\documentclass[11pt]{amsart}
\usepackage{amsfonts}
\usepackage{}
\usepackage{amsfonts}
\usepackage[arrow,matrix]{xy}
\usepackage{amsmath,amssymb,bbm, amscd, amsthm,mathrsfs,hyperref}
\usepackage{longtable}
\theoremstyle{plain}
\newtheorem{thm}{Theorem}[section]
\newtheorem{cor}[thm]{Corollary}
\newtheorem{lem}[thm]{Lemma}
\newtheorem{prop}[thm]{Proposition}
\theoremstyle{definition}
\newtheorem{defn}[thm]{Definition}
\newtheorem{eg}[thm]{Example}

\newtheorem{rem}[thm]{Remark}

\def\al{\alpha}
\def\bt{\beta}
\def\dt{\delta}

\def\sg{\sigma}
\def\tt{\theta}
\def\vph{\varphi}
\def\lmd{\lambda}
\def\vps{\varepsilon}
\def\ome{\omega}

\def\bgm{\Gamma}

\def\ms{\mathscr}
\def\mc{\mathcal}
\def\mb{\mathbb}
\def\mk{\mathfrak}

\def\ot{\otimes}

\def\se{\leqslant}
\def\le{\geqslant}
\def\lan{\langle}
\def\ran{\rangle}

\def\Hom{\operatorname {Hom}}

\def\RHom{\operatorname {RHom}}
\def\Ext{\operatorname {Ext}}
\def\Tor{\operatorname {Tor}}

\def\dim{\operatorname {dim}}

\def\id{\operatorname {id}}

\def\HH{\operatorname {HH}}

\def\injdim{\operatorname {injdim}}
\def\projdim{\operatorname {projdim}}
\def\Gr{{\operatorname {Gr}}}
\def\ad{{\operatorname {ad}}}

\def\t{\text}

\def\it{\textit}

\def\iff{if and only if}
\def\kk{\mathbbm{k}}

\def\cy{Calabi-Yau}

\def\ra{\rightarrow}
\def\xra{\xrightarrow}

\def\ZZ{\mathbb{Z}}
\def\NN{\mathbb{N}}
\def\rr{\mathbbm{r}}

\begin{document}
\title{Calabi-Yau pointed Hopf algebras of finite Cartan type}

\author{Xiaolan YU}
\address {Xiaolan YU\newline Hangzhou Normal University, Hangzhou, Zhejiang 310036, China}
\email{xlyu1005@yahoo.com.cn}

\author{Yinhuo ZHANG}
\address {Yinhuo ZHANG\newline Department WNI, University of Hasselt, Universitaire Campus, 3590 Diepeenbeek,Belgium } \email{yinhuo.zhang@uhasselt.be}

\date{}

\begin{abstract}
We study the Calabi-Yau property of pointed Hopf algebra $U(\mc{D},\lmd)$ of finite Cartan type. It turns out that this class of pointed Hopf algebras constructed by N. Andruskiewitsch and H.-J. Schneider contains many Calabi-Yau Hopf algebras. To give concrete  examples of new Calabi-Yau Hopf algebras, we classify the  Calabi-Yau pointed Hopf algebras $U(\mc{D},\lmd)$ of dimension less than 5.
\end{abstract}

\maketitle

\section*{Introduction}

In \cite{as3},  N. Andruskiewitsch and H.-J. Schneider classified
pointed Hopf algebras with finite Gelfand-Kirillov dimension, which are
domains with a finitely generated abelian group of group-like elements, and
a positive braiding.
In the same paper, the authors constructed a class of Hopf algebra $U(\mc{D},\lmd)$,
generalizing the quantized enveloping algebra $U_q(\mathfrak{g})$ of
a finite dimensional semisimple Lie algebra $\mathfrak{g}$. These pointed Hopf algebras turn out to be Artin-Schelter (AS-) Gorenstein Hopf algebras. AS-Gorenstein Hopf algebras have been recently intensively studied (e.g. \cite{bro,bz,lwz,lwz1,wz,wz1}). One of important
properties of an AS-Gorenstein Hopf algebra is the existence of a homological integral, which generalizes Sweedler's classical integral of a finite dimensional Hopf algebra cf.\cite{lwz}.  Brown and Zhang proved that the rigid dualizing complex of an AS-Gorenstein  Hopf algebra is determined by its Homological integral and antipode \cite{bz}. He, Van Oystaeyen and Zhang used homological integrals to investigate the Calabi-Yau property of cocommutative Hopf algebras \cite{hvz}. They  successfully classified the low dimensional cocommutative Calabi-Yau Hopf algebras over an algebraic closed field of characteristic zero.

The main aim of this paper is to find out when a pointed Hopf algebra
$U(\mc{D},\lmd)$ of finite Cartan type is  Calabi-Yau, and to classify those low dimensional Calabi-Yau pointed Hopf algebras. It turns out that the class $U(\mc{D},\lmd)$ of pointed Hopf algebras contains many Calabi-Yau Hopf algebras. Most of them are of different types from the quantum groups $U_q(\mathfrak{g})$, which were proved to be Calabi-Yau by Chemla \cite{c}.  This give us more interesting examples of Calabi-Yau Hopf algebras. The paper is organized as follows.

In Section 1, we recall the pointed Hopf algebras $U(\mc{D},\lmd)$ of finite Cartan type, the definition of a Calabi-Yau algebra, the notion of a homological integral and a (rigid) dualizing complex over a Noetherian algebra. In Section 2, we study the Calabi-Yau pointed Hopf algebras of finite Cartan type.
We give a necessary and sufficient condition for a Hopf algebra  $U(\mc{D},\lmd)$ to be Calabi-Yau,
and calculate the rigid dualizing complex of $U(\mc{D},\lmd)$ (Theorem \ref{cy}). In \cite{c}, the author calculated the rigid dualizing complexes of quantum groups $U_q(\mathfrak{g})$. As a
consequence of the characterization theorem, the quantum groups $U_q(\mathfrak{g})$ are Calabi-Yau Hopf algebras.

 A pointed Hopf algebra of form $U(\mc{D},\lmd)$ is a (cocycle) deformation of the smash product $\mc{B}(V)\#\kk\bgm$, where $\mc{B}(V)$ is the Nichols algebra of $V$ and $\kk \bgm$ is the group algebra of the group formed by group-like elements of $U(\mc{D},\lmd)$. Our second aim in this paper is to study the Calabi-Yau property of the Nichols algebra $\mc{B}(V)$. The algebra $\mc{B}(V)$ is an $\NN^{p+1}$-filtered algebra. By analyzing the rigid dualizing complex of the associated graded algebra $\mathbb{G}\rr \mc{B}(V)$, we obtain the rigid dualizing complex of  $\mc{B}(V)$. We then are able to give a necessary and sufficient condition for the algebra $\mc{B}(V)$ to be Calabi-Yau, which forms the main result of Section 3 (see Theorem \ref{nichols cy}).

In Section 4, we discuss the relation between the Calabi-Yau property of a pointed Hopf algebra $U(\mc{D},\lmd)$ and the Calabi-Yau property of the associated Nichols algebra $\mc{B}(V)$. It turns out that for a pointed Hopf algebra $U(\mc{D},\lmd)$ and its associated Nichols algebra $\mc{B}(V)$, if one of them is CY, the other one is not.

In the final section, we classify the
Calabi-Yau pointed Hopf algebras $U(\mc{D},\lmd)$ of dimension less
than 5. It turns out that $U_q(\mathfrak{sl}_2)$ is the only known  non-cocommutative example in the classification. The other non-cocommutative Hopf algebras are new examples of Calabi-Yau Hopf algebras.

\section{Preliminaries}
Throughout this paper, we fix an algebraically closed field $\kk$. All
vector spaces, algebras are over $\kk$. The unadorned  tensor $\ot$
means $\ot_\kk$. Given an algebra $A$, we write $A^{op}$ for the
opposite algebra of $A$ and $A^e$ for the enveloping algebra $A\ot
A^{op}$ of $A$.

Let $A$ be an algebra. For a left $A$-module $M$ and an algebra automorphism $\phi:A\ra A$,
$_{\phi}M$ stands for the left $A$-module twisted by the
automorphism $\phi$. Similarly, for a right $A$-module $N$, we have
$N_{\phi}$. Observe that  $A_\phi\cong {}_{\phi^{-1}} A$ as
$A$-$A$-bimodules. $A_\phi\cong A$ as $A$-$A$-bimodules if and only
if $\phi$ is an inner automorphism.

 Let $A$ be a Hopf algebra, and $\xi:A\ra\kk$ an algebra
homomorphism. We write $[\xi]$ to be the \textit{winding} homomorphism of
$\xi$ defined by
$$[\xi](a)=\sum\xi(a_1)a_2,$$ for any $a\in A$.

A Noetherian algebra in this paper means a \textit{left and right} Noetherian algebra.

\subsection{Pointed Hopf algebra $U(\mc{D},\lmd)$} In this subsection we recall the definitions and basic properties of
Nichols algebras and pointed Hopf algebras of finite Cartan type. More details can be found in \cite{as3}.
We fix the following terminology.
\begin{itemize}
\item  a free abelian group $\bgm$ of finite rank $s$;
\item  a Cartan matrix $(a_{ij})\in \mathbb{Z}^{\tt\times \tt}$ of finite type, where $\tt\in\NN$. Denote by $(d_1,\cdots, d_\tt)$ a diagonal matrix
of positive integers such that $d_ia_{ij} = d_ja_{ji}$, which is
minimal with this property;
\item a set $\mathcal {X}$ of connected components of the Dynkin diagram corresponding
to the Cartan matrix $(a_{ij})$. If $1\se i, j\se \tt$, then $i\sim
j$ means that they belong to the same connected component;
\item a family $(q_{_I})_{I\in \mc{X}}$ of elements in $\kk$ which are \textit{not} roots of unity;
\item elements $g_1,\cdots , g_\tt\in \bgm$ and characters $\chi_1,\cdots, \chi_\tt\in \hat{\bgm}$ such that
\begin{equation}\label{q}\chi_j(g_i)\chi_i(g_j)=q_{I}^{d_ia_{ij}}, \t{  } \chi_i(g_i)=q_{I}^{d_i}, \t{   for all $1\se i,j\se \tt$, $I\in
\mc{X}$}.\end{equation}
\end{itemize}

Let $\mc{D}$ be the collection $\mc{D}(\bgm,(a_{ij})_{1\se i,j\se
\tt}, (q_{_{I}})_{I\in \mathcal {X}}, (g_i)_{1\se i\se
\tt},(\chi_i)_{1\se i\se \tt} )$. A \textit{linking datum}
$\lmd=(\lmd_{ij})$ for $\mc{D}$ is a collection of elements
$(\lmd_{ij})_{1\se i<j\se \tt,i\nsim j}\in \{0,1\}$ such that
$\lmd_{ij}=0$ if $g_ig_j=1$ or $\chi_i\chi_j\neq\varepsilon$. We
write the datum $\lmd=0$, if $\lmd_{ij}=0$ for all $1\se i<j\se
\tt$. The datum $(\mc{D},\lmd)=(\bgm,(a_{ij}) , q_{_I},
(g_i),(\chi_i), (\lmd_{ij}) )$ is called a \textit{ generic datum of
finite Cartan type} for group $\bgm$.

\begin{defn}\cite[Sec. 4]{as3}\label{ud} Let $(\mc{D},\lmd)$ be a generic datum of finite Cartan type.
Let  $U(\mc{D},\lmd)$ be the algebra generated by $x_1,\cdots,x_\tt$
and $y_1^{\pm1},\cdots,y_s^{\pm1}$ subject to the relations
$$\begin{array}{lcr}
&y_m^{\pm1}y_h^{\pm1}=y_h^{\pm1}y_m^{\pm1},\;\;\;y_m^{\pm1}y_m^{\mp1}=1,&1\se m,h\se s,\\
\t{(group action)}&y_hx_j=\chi_j(y_h)x_jy_h,&1\se j\se \tt,\;\;1\se h\se s,\\
\t{(Serre relations)}&(\t{ad}_cx_i)^{1-a_{ij}}(x_j)=0,&1\se i\neq j\se \tt,\;\;i\sim j,\\
\t{(linking
relations)}&x_ix_j-\chi_j(g_i)x_jx_i=\lmd_{ij}(1-g_ig_j),&1\se
i<j\se \tt,\;\;i\nsim j,
\end{array}$$
where $\t{ad}_c$ is the braided adjoint representation defined  in
\cite[Sec. 1]{as3}.
\end{defn}

For a generic datum of finite Cartan type $(\mc{D},\lmd)$,   denote
by $q_{ji}=\chi_i(g_j)$. Then Equation (\ref{q}) reads as follows:
\begin{equation}\label{q1}q_{ii}=q_I^{d_i}\t{ and } q_{ij}q_{ji}=q_{I}^{d_ia_{ij}}\t{ for all }
1\se i,j\se \tt, I\in \mc{X}.
\end{equation}
Let $V$ be a Yetter-Drinfeld module
over the group algebra $\kk \bgm$ with basis $x_i\in
V^{\chi_i}_{g_i}$, $1\se i\se \tt$. In other words, $V$ is a braided
vector space with basis $x_1,\cdots, x_\tt$, whose braiding is given
by
$$c(x_i\ot x_j)=q_{ij}x_j\ot x_i.$$ It can be easily derived from the proof of \cite[Thm. 4.3]{as3} that the Nichols algebra $\mathcal {B}(V)$ is isomorphic to the following algebra:
$$ \kk\langle x_1, \cdots , x_\tt\mid (\t{ad}_cx_i)^{1-a_{ij}}(x_j)= 0, 1 \se i, j \se \tt, i\neq
j\rangle.$$
We refer to \cite[Sec. 2]{as} for the definition of a
Nichols algebra.

Let $\Phi$ be the root system corresponding to the Cartan matrix
$(a_{ij})$ with $\{\al_1,\cdots, \al_\tt\}$ a set of fix simple
roots, and  $\mathcal {W}$ the Weyl group. We fix a reduced
decomposition of the longest element $w_0=s_{i_1}\cdots s_{i_p}$ of
$\mathcal {W}$ in terms of the simple reflections.  Then the
positive roots are precisely the followings,
$$\bt_1=\al_{i_1}, \;\;\bt_2=s_{i_1}(\al_{i_2}),\cdots, \bt_p=s_{i_1}\cdots s_{i_{p-1}}(\al_{i_p}).$$ If
$\bt_i=\sum_{i=1}^{\tt} m_i\al_i$, then we write
$$g_{\bt_i}=g_1^{m_1}\cdots g_\tt^{m_\tt} \t{ and }\chi_{\bt_i}={\chi}_1^{m_1}\cdots {\chi}_\tt^{m_\tt}.$$ Similarly,
we write $q_{_{\bt_j\bt_i}}=\chi_{_{\bt_i}}(g_{_{\bt_j}})$.

Root vectors for a quantum group $U_q(\mathfrak{g})$ were defined by
Lusztig \cite{l}. Up to a non-zero scalar, each root vector can be
expressed as an iterated braided commutator. As in \cite[Sec.
4.1]{as4}, this definition can be generalized to a pointed Hopf
algebras $U(\mc{D},\lmd)$. For each positive root $\bt_i$, $1\se i
\se p$, the root vector $x_{\bt_i}$ is defined by the same iterated
braided commutator of the elements $x_1, \cdots , x_\tt$, but with
respect to the general braiding.

\begin{rem}\label{root} If $\bt_j=\al_l$, then we have $x_{\bt_j}=x_l$. That is, $x_1,\cdots,x_\tt$ are the simple root vectors.
\end{rem}
\begin{lem}\cite[Thm. 4.3]{as3}\label{pointeddef}
Let $(\mc{D},\lmd)=(\bgm,(a_{ij}) , q_{_I}, (g_i),(\chi_i),
(\lmd_{ij}) )$ be a generic datum of finite Cartan type for $\bgm$.
The algebra $U(\mc{D},\lmd)$ as defined in Definition \ref{ud} is a
pointed Hopf algebra with comultiplication structure determined by
$$\Delta(y_h)=y_h\ot y_h,\;\; \Delta (x_i)=x_i\ot 1+g_i\ot x_i,\;\;1\se h\se s,1\se i\se \tt.$$
Furthermore, $U(\mc{D},\lmd)$ has a PBW-basis given by monomials in
the root vectors
\begin{equation}\label{pbw}\{x_{\bt_1}^{a_1}\cdots
x_{\bt_p}^{a_p}y\},\end{equation} for $a_i\le 0$,  $1\se i\se p$,
and $y\in \bgm$. The coradical filtration of $U(\mc{D},\lmd)$ is
given by
$$U(\mc{D},\lmd)_N = \t{span of }x_{\bt_{i_1}}\cdots x_{\bt_{i_j}}y,\;\;j\se N,\;\;y\in  \bgm.
$$There is an isomorphism of graded Hopf algebras $\Gr U(\mc{D}, \lmd)\cong \mc{B}(V)\#\kk \bgm\cong U(\mc{D}, 0)$.  The algebra $U(\mc{D},\lmd)$ has finite Gelfand-Kirillov
dimension and is a domain.

\end{lem}

In \cite{as3}, degrees of the PBW basis elements are defined as follows:
\begin{equation}\label{deg pointd}\deg (x_{_{\bt_1}}^{a_1}\cdots x_{_{\bt_p}}^{a_p}y)=(a_1,\cdots,
a_p,\sum_{i=1}^p a_iht(\bt_i))\in
(\mathbb{Z}^{\le0})^{p+1},\end{equation}
where $ht(\bt)$ is the height of the root $\bt$. That is, if $\bt=\sum_{i=1}^\tt
m_i\al_i$, then $ht(\bt)=\sum_{i=1}^\tt m_i$. Order the elements
in $(\mathbb{Z}^{\le0})^{p+1}$ as follows
\begin{equation}\label{ordering}\begin{array}{l}
(a_1, \cdots, a_p, a_{p+1})<(b_1, \cdots, b_p, b_{p+1})\t{ if and
only if there is some } \\1\se k\se p+1 ,\t{ such that } a_i=b_i \t{
for }i\le k \t{ and }a_{k-1}<b_{k-1}.\end{array}
\end{equation}
Given $m\in \NN^{p+1}$, let $F_m U(\mc{D},\lmd)$ be the space spanned
by the monomials $x_{\bt_1}^{a_1}\cdots x_{\bt_p}^{a_p}y$ such that
$\deg(x_{\bt_1}^{a_1}\cdots x_{\bt_p}^{a_p}y)\se m$. Then we obtain a filtration on the algebra $U(\mc{D},\lmd)$.

\begin{lem}\label{filt1}
If the root vectors $x_{\bt_i}, x_{\bt_j}$ belong to the same
connected component and  $j>i$, then
\begin{equation}\label{comm}[x_{\bt_i}, x_{\bt_j}]_c=\sum_{\textbf{a}\in \mathbb{N}^p}\rho_{\textbf{a}}x_{\bt_1}^{a_1}\cdots x_{\bt_p}^{a_p}, \end{equation}
where $\rho_{\textbf{a}}\in \kk$ and $\rho_{\textbf{a}}\neq0$ only
when $\textbf{a}=(a_1,\cdots,a_p)$ is such that $a_s=0$ for $s\se i$
or $s\le j$. In particular, in $U(\mc{D},0)$, the equation
(\ref{comm}) holds for all root vectors $x_{\bt_i}, x_{\bt_j}$ with
$i<j$.
\end{lem}
\proof This follows from \cite[Prop. 2.2]{as3} and the classical
relations that hold for a quantum group $U_q(\mathfrak{g})$ (see
\cite[Thm. 9.3]{dp} for example). It was actually proved in  Step VI
of the proof of Theorem 4.3 in \cite{as3}.  \qed

\begin{lem}\label{filt}
The filtration defined by PBW basis is an algebra filtration. The
associated graded algebra $\mathbb{G}\mathbbm{r}U(\mc{D},\lmd)$ is
generated by $x_{\bt_i}$, $1\se i\se p$, and $y_h$, $1\se h\se s$,
subject to the relations
$$y_h^{\pm1}y_m^{\pm1}=y_m^{\pm1}y_h^{\pm1},\;\;\;y_h^{\pm1}y_h^{\mp1}=1,\;\;\; 1\se h,m\se s,$$
$$y_hx_{\bt_i}=\chi_{\bt_i}(y_h)x_{\bt_i} y_h,\;\;\;1\se i\se p,\;\;1\se h\se s,$$
$$x_{\bt_i}x_{\bt_j}=\chi_{_{\bt_j}}(g_{_{\bt_i}})x_{\bt_j}x_{\bt_i},\;\;1\se
i<j\se p.$$
\end{lem}
\proof This follows from Lemma \ref{filt1} and the linking
relations. \qed

\subsection{Calabi-Yau algebras}Following \cite{g2},  we call an algebra  \textit{\cy\ of dimension
$d$} if
\begin{enumerate}\item[(i)] $A$ is homologically smooth, that is, $A$ has
a bounded resolution of finitely generated projective
$A$-$A$-bimodules;
\item[(ii)] There are $A$-$A$-bimodule
isomorphisms$$\Ext_{A^e}^i(A,A^e)=\begin{cases}0& i\neq d;
\\A&i=d.\end{cases}$$
\end{enumerate}
In the sequel, Calabi-Yau will be abbreviated to CY for short.

In \cite{hvz}, the CY property of Hopf algebras was discussed by
using the homological integrals of Artin-Schelter Gorenstein (AS-Gorenstein for short)
algebras \cite[Thm. 2.3]{hvz}.

 Let us recall the definition of  an AS-Gorenstein algebra
(cf. \cite{bz}). \label{defn as}\index{AS-regular}
\begin{enumerate}
\item[(i)] Let $A$ be a  left Noetherian augmented algebra with
a fixed augmentation map $\vps:A\ra \kk$. The algebra $A$ is said to be \it{left
AS-Gorenstein},\index{AS-Gorenstein!left} if
\begin{enumerate}
\item[(a)] $\injdim  {_AA}=d<\infty$,
\item[(b)] $\dim\Ext_A^i({_A\kk},{_AA})=\begin{cases}0,&i\neq d;\\ 1,& i=d,\end{cases}$
\end{enumerate}
where $\injdim$ stands for  injective dimension.

A \textit{Right AS-Gorenstein algebras}\index{AS-Gorenstein!right} can be defined similarly.

\item[(ii)] An algebra $A$ is
said to be \textit{AS-Gorenstein}\index{AS-Gorenstein} if it is both left and right
AS-Gorenstein (relative to the same augmentation map $\vps$).
\item[(iii)] An
AS-Gorenstein algebra $A$ is said to be \it{regular}  if, in
addition, the global dimension of $A$
 is finite.
\end{enumerate}

Let $A$ be a Noetherian algebra. If the injective dimension of $_AA$ and $A_A$
are both finite, then these two integers are equal by \cite[Lemma
A]{za}. We call this common value the injective dimension of $A$.
The left global dimension and the right global dimension of a
Noetherian algebra are equal \cite[Exe. 4.1.1]{we}.  When the global
dimension is finite, then it is equal to the injective dimension.

In order to study infinite dimensional Noetherian Hopf algebras, Lu,
Wu and Zhang introduced the concept of a homological integral for an
AS-Gorenstein Hopf algebra in \cite{lwz}, which is a generalization
of an integral of a finite dimensional Hopf algebra. In \cite{bz},
homological integrals were defined for general AS-Gorenstein
algebras.

 Let $A$ be a left AS-Gorenstein algebra
with $\injdim{_AA}=d$. Then $\Ext_A^d({_A\kk},{_AA})$ is a 1-dimensional
right $A$-module. Any nonzero element in $\Ext_A^d({_A\kk},{_AA})$
is called a \it{left homological integral} of $A$. We write
$\int_A^l$ for $\Ext_A^d({_A\kk},{_AA})$. Similarly, if $A$ is right
AS-Gorenstein, any nonzero element  in $\Ext_A^d({\kk_A},{A_A})$ is
called a \it{right homological integral} of $A$. Write
$\int_A^r$ for $\Ext_A^d({\kk_A},{A_A})$.

$\int_A^l$ and $\int_A^r$
are called\textit{ left and right homological integral modules} of $A$
respectively.

CY algebras are closely related to algebras having rigid dualizing complexes. The
non-commutative version of a dualizing complex was first introduced
by Yekutieli.

\begin{defn}\cite{y} (cf. \cite[Defn. 6.1]{vdb})\label{defn dc} \index{dualizing complex}
Assume that $A$ is a  (graded) Noetherian algebra. Then an object
$\ms{R}$ of $D^b(A^e)$ ($D^b(\t{GrMod}(A^e))$) is called a
\textit{dualizing complex} (in the graded sense) if it satisfies the
following conditions:
\begin{enumerate}
\item[(i)] $\ms{R}$ is of  finite injective dimension over $A$ and $A^{op}$.
\item[(ii)] The cohomology of $\ms{R}$ is given by bimodules which are
finitely generated on both sides.
\item[(iii)] The natural morphisms $A\ra  {\RHom}_A(\ms{R},\ms{R})$ and $A\ra
\RHom_{A^{op}}(\ms{R},\ms{R})$ are isomorphisms in $D(A^e)$
($D(\t{GrMod}(A^e))$).
\end{enumerate}
\end{defn}

Roughly speaking, a dualizing complex is a complex $\ms{R}\in
D^b(A^e)$ such that the functor
\begin{equation}\label{dualcom}\RHom_A(-,\ms{R}):D^b_{fg}(A)\ra D^b_{fg}(A^{op})\end{equation} is a
duality, with adjoint $\RHom_{A^{op}}(-,\ms{R})$ (cf. \cite[Prop.
3.4 and Prop. 3.5]{y}). Here $D^b_{fg}(A)$ is the  full triangulated
subcategory of $D(A)$ consisting of bounded complexes with finitely
generated cohomology modules.

In the above definition,  the algebra $A$ is a Noetherian algebra.
In this case, a dualizing complex in the graded sense is also a
dualizing complex in the usual sense.

Dualizing complexes are not unique  up to isomorphism. To overcome
this weakness,  Van den Bergh introduced the concept of a rigid
dualizing complex in \cite[Defn. 8.1]{vdb}.

\begin{defn}\label{defn rdc}\index{dualizing complex!rigid} Let $A$ be a (graded) Noetherian algebra. A dualizing
 complex $\ms{R}$ over $A$ is called \textit{rigid} (in the graded sense) if $$\RHom_{A^e}(A,{_A\ms{R}\ot \ms{R}_A})\cong
 \ms{R}$$ in $D(A^e)$ ($D(\t{GrMod}(A^e))$).
\end{defn}

Note again that if $A^e$ is Noetherian then the graded version of
this definition implies the ungraded version.

\begin{lem}\label{vdb} \rm{(cf. \cite[Prop. \;4.3]{bz}  and  \cite[Prop. \;8.4]{vdb})} Let $A$ be a Noetherian algebra.   Then the following two conditions are
equivalent:
\begin{enumerate}\item[(a)] $A$ has a rigid dualizing complex $\ms{R}=A_{\psi}[s]$, where $\psi$ is an
algebra automorphism and $s\in \ZZ$.
\item[(b)] $A$ has finite
injective dimension $d$ and there is an algebra automorphism $\phi$
such that
$$\Ext_{A^e}^i(A,A^e)\cong\begin{cases}0,& i\neq d;
\\A_\phi&i=d\end{cases}$$ as $A$-$A$-bimodules. \end{enumerate}
In this case,  $\phi=\psi^{-1}$ and $s=d$.
\end{lem}

 The
following corollary follows immediately from Lemma \ref{vdb} and the
definition of a CY algebra. It characterizes the Noetherian CY algebras.

\begin{cor}\label{cor cyrid}
Let $A$ be a Noetherian algebra which is  homologically smooth. Then
$A$ is a CY algebra of dimension $d$ if and only if $A$ has a rigid
dualizing complex $A[d]$.
\end{cor}

\section{Calabi-Yau pointed Hopf algebras of finite Cartan type}

In this section, we  calculate the rigid dualizing complex of a
pointed Hopf algebra $U(\mc{D},\lmd)$ and  study its Calabi-Yau
property.

Before we give the main theorem of this section, let us recall the Koszul complex of quadratic algebras (cf. \cite{vdb1}). Let $V$ be a finite dimensional vector space and $T(V)$ the tensor algebra of $V$. Suppose that $A$ is a quadratic algebra, that is, $A=T(V)/\lan R \ran$, where $R\subseteq V\ot V$. The quadratic dual algebra of $A$, denoted by $A^!$, is the quadratic algebra $T(V^*)/\lan R^ \perp\ran$. Let $\{x_i\}_{i=1,\cdots,n}$ be a basis of $V$ and $\{x^*_i\}_{i=1,\cdots,n}$ be the dual basis of $V^*$. Introduce the canonical element $e=\sum_{i=1}^n x_i\ot x^*_i\in A\ot A^!$. The right multiplication by $e$ defines a complex
\begin{equation}\label{koszulcom} \cdots \ra A\ot
A^{!*}_{j}\xra{d_j} A\ot A^{!*}_{j-1}\ra\cdots\ra A\ot A^{!*}_{1}\ra
A\ra \kk\ra 0.\end{equation}

This complex is called the Koszul complex of $A$.
The algebra $A$ is Koszul if and only if the complex (\ref{koszulcom}) is a resolution of $_A\kk$.

Let $\mathcal {K}$ be the bimodule complex defined as follows
\begin{equation}\label{kosbimod}\mathcal {K}: \cdots \ra A\ot
A^{!*}_{j}\ot A\xra{D_j} A\ot A^{!*}_{j-1}\ot A\ra \cdots \ra A\ot
A\ra 0.\end{equation} The differentials $D_j:A\ot A^{!*}_{j}\ot A\ra
A\ot A^{!*}_{j-1}\ot A$, $1\se j\se n$,  are defined by
$D_j=d_j^l+(-1)^jd_j^r$, where $d_j^l(1\ot a\ot1)=\sum_{i=1}^n
x_i\ot  a\cdot x^*_i\ot1$ and $d_j^r(1\ot a\ot1)=\sum_{i=1}^n 1\ot
x^*_i\cdot a \ot x_i$, for any $1\ot a\ot 1\in A\ot A^{!*}_{j}\ot
A$. The complex $\mathcal {K}$ is called the Koszul bimodule complex of $A$. If $A$ is Koszul, then  $\mathcal {K}\ra A\ra 0$ is exact.

In the rest of this section we fix a generic  datum of finite Cartan
type  $$(\mc{D},\lmd)=(\bgm,(a_{ij})_{1\se i,j\se \tt},(q_{_I})_{I\in \mathcal {X}},
(g_i)_{1\se i\se \tt},(\chi_i)_{1\se i\se \tt},(\lmd_{ij})_{1\se
i<j\se \tt,i\nsim j} ),$$ where $\bgm$  is a free abelian group of
rank $s$. Let $x_{\bt_1},\cdots, x_{\bt_p}$ be the root vectors.
Recall from Remark \ref{root} that there are $1\se j_k\se p$, $1\se k\se
\tt$, such that $x_{\bt_{j_k}}=x_k$.

\begin{lem}\label{global}
The Hopf algebra $A=U(\mc{D},\lmd)$ is Noetherian with finite global
dimension bounded by $p+s$.
\end{lem}
\proof  The group algebra $\kk \bgm$ is isomorphic to a Laurent
polynomial algebra with $s$ variables. So $\kk \bgm$ is Noetherian
of global dimension $s$. By Lemma \ref{filt1}, the algebra $\Gr A\cong
U(\mc{D},0)$ is an iterated Ore extension of $\kk \bgm$. Indeed, if
$x_{\bt_1},\cdots,x_{\bt_p}$ are the root vectors of $A$, then
$$\Gr A\cong
\kk\bgm[x_{\bt_1};\tau_1,\delta_1][x_{\bt_2};\tau_2,\delta_2]\cdots[x_{\bt_p};\tau_p,\delta_p],$$
where for $1\se j\se p$, $\tau_j$ is an algebra automorphism such
that $\tau_{j}(x_{\bt_i})$ is just a scalar multiple of $x_{\bt_i}$
for $i<j$, and $\delta_j$ is a $\tau_{j}$-derivation such that
$\delta_j(x_{\bt_i})$, $i<j$, is a linear combination of monomials
in $x_{\bt_{i+1}},\cdots, x_{\bt_{j-1}}$. By \cite[Thm. 1.2.9 and
Thm. 7.5.3]{mr},  we have that $\Gr A$ is Noetherian of global
dimension less than $p+s$. Now it follows from  \cite[Thm. 1.6.9 and
Cor. 7.6.18]{mr} that the algebra $A$ is Noetherian of global
dimension less than $p+s$. \qed

\begin{thm}\label{as}
Let ($\mc{D}, \lmd)$ be a generic datum of finite Cartan  and  $A$
the Hopf algebra $U(\mc{D},\lmd)$. Then $A$ is Noetherian AS-regular
of global dimension $p+s$, where $s$ is the rank of $\bgm$ and $p$
is the number of the positive roots of the Cartan matrix. The left
homological integral module $\int^l_A$ of $A$ is isomorphic to
$\kk_\xi$, where $\xi:A\ra \kk$ is an algebra homomorphism defined
by $\xi(g)=(\prod_{i=1}^p\chi_{_{\bt_i}})(g)$ for all $g\in \bgm$
and $\xi(x_i)=0$ for all $1\se i\se \tt$.
\end{thm}
\proof  We first show that $$\Ext_{A}^i(_A\kk,{_AA})\cong\begin{cases}0,&i\neq p+s;\\
\kk_\xi ,&i=p+s.\end{cases}$$ With Lemma \ref{filt} and Lemma
\ref{global}, the method in \cite[Prop. 3.2.1]{c} for computing the
group
$\Ext_{U_q(\mathfrak{g})}^*(_{U_q(\mathfrak{g})}\kk,{_{U_q(\mathfrak{g})}U_q(\mathfrak{g})})$
also works in the case of $A=U(\mc{D},\lmd)$. The  difference is that
the right $A$-module structure on $\Ext_{A}^{p+s}(_A\kk,{_AA})$ is
not trivial in the case of $U(\mc{D},\lmd)$. Let
$C=\mathbb{G}\mathbbm{r}U(\mc{D},\lmd)$. We also have that
$\Ext_{A}^i(_A\kk,{_AA})=0$ for $i\neq p+s$ and
$\Ext_{C}^{p+s}(_C\kk,{_CC})\cong\Ext_{A}^{p+s}(_A\kk,{_AA})$ as
right $\bgm$-modules.

We now  give the structure of $\Ext_{C}^*(_C\kk,{_CC})$. Let $B$ be
the following algebra,
$$\kk\langle x_{_{\bt_1}}, \cdots,x_{_{\bt_p}}\mid x_{_{\bt_i}}x_{_{\bt_j}}=\chi_{_{\bt_j}}(g_{_{\bt_i}})x_{_{\bt_j}}x_{_{\bt_i}},
1\se i<j\se p\rangle.$$ Then $C=B\#\kk\bgm$. We have the following
isomorphisms
$$\begin{array}{ccl}\RHom_C(\kk,C)&\cong&\RHom_C(\kk\bgm\ot_{\kk\bgm}\kk,C)\\
&\cong&\RHom_{\kk\bgm}(\kk,\RHom_C(\kk\bgm,C))\\
&\cong&\RHom_{\kk\bgm}(\kk,\kk\bgm)\ot^L_{\kk\bgm}\RHom_C(\kk\bgm,C).
\end{array}$$

Let
\begin{equation}\label{koszulcom1}0\ra B\ot B^{!*}_p \ra \cdots \ra B\ot B^{!*}_{j} \ra
\cdots\ra B\ot B^{!*}_{1} \ra B\ \ra \kk \ra 0 \end{equation} be
 the Koszul complex of $B$ (cf. complex (\ref{koszulcom})). It is a
 projective resolution of $\kk$.
Each $B_j^{!*}$
is a left $\kk\bgm$-module defined by
$$\begin{array}{ccl}[g(\bt)](x^*_{\bt_{i_1}}\wedge\cdots\wedge
x^*_{\bt_{i_j}})&=&\bt(g^{-1}(x^*_{\bt_{i_1}}\wedge\cdots\wedge
x^*_{\bt_{i_j}}))\\&=&\bt(g^{-1}(x^*_{\bt_{i_1}})\wedge\cdots\wedge
g^{-1}(x^*_{\bt_{i_j}}))\\
&=&\prod_{t=1}^j\chi_{\bt_{i_t}}(g)\bt(x^*_{\bt_{i_1}}\wedge\cdots\wedge
x^*_{\bt_{i_j}}).\end{array}$$ Thus, each $B\ot B_j^{!*}$ is a
$B\#\kk\bgm$-module defined by $$(c\#g)\cdot(b\ot\bt )=(c\#g)(b)\ot
g(\bt),$$ for any $b\ot \bt\in B\ot B^{!*}_j$ and $c\#g\in
B\#\kk\bgm$. It is not difficult to see that the complex
(\ref{koszulcom1}) is an exact sequence of $B\#\kk\bgm$ modules.
 Tensoring  it  with $\kk\bgm$, we obtain the following exact sequence of  $B\#\kk\bgm$-modules
$$0\ra B\ot B^{!*}_p\ot\kk\bgm\ra \cdots \ra B\ot
B^{!*}_{j}\ot\kk\bgm\ra \cdots\ra B\ot B^{!*}_{1}\ot\kk\bgm\ra
B\ot\kk\bgm\ra \kk\bgm\ra 0,$$ where the $\bgm$-action is diagonal.
Each $B\ot B^{!*}_j\ot\kk\bgm$ is a free $B\#\kk\bgm$-module.
Therefore, we obtain a projective resolution of $\kk\bgm$ over
$B\#\kk\bgm$.

The complex
$$0\ra \Hom_{C}(B\ot\kk\bgm,C)\ra \Hom_C(B\ot
B^{!*}_{1}\ot\kk\bgm,C)\ra \cdots \ra\Hom_C(B\ot
B^{!*}_p\ot\kk\bgm,C) \ra 0$$ is isomorphic to the following complex
$$0\ra C\ra
B^{!}_1\ot C\ra \cdots \ra B^{!}_{p-1}\ot C\xra{\dt_p} B^{!}_{p}\ot
C\ra 0.$$ This complex is exact except at $B^{!}_{p}\ot C$, whose
cohomology is isomorphic to $B^{!}_{p}\ot \kk\bgm$. So
$\RHom_C(\kk\bgm,C)\cong B_p^!\ot\kk\bgm[p]$. We have
$$(x^*_{_{\bt_1}}\wedge\cdots\wedge x^*_{_{\bt_p}})\ot
g=(\prod_{i=1}^p\chi_{_{\bt_i}})(g)g((x^*_{_{\bt_1}}\wedge\cdots\wedge
x^*_{_{\bt_p}})\ot 1),$$ for all $g\in \bgm$. The group $\bgm$ is a free abelian group of rank $s$, so
$\RHom_{\kk\bgm}(\kk,\kk\bgm)\cong \kk[s]$. Therefore,  we obtain
that $$\RHom_{\kk\bgm}(\kk,\kk\bgm)\ot^L_{\kk\bgm}\RHom_C(\kk\bgm,C)\cong \kk_{\xi'}[p+s],$$
where $\xi'$ is defined by
$\xi'(g)=(\prod_{i=1}^p\chi_{_{\bt_i}})(g)$ for all $g\in \bgm$ and
$\xi'(x_{_{\bt_j}})=0$ for all $1\se j\se p$. That is,
$$\Ext_{C}^i(_C\kk,{_CC})\cong\begin{cases}0,&i\neq p+s;\\
\kk_{\xi'} ,&i=p+s.\end{cases}$$

$\Ext_{A}^{p+s}(_A\kk,{_AA})$ is a 1-dimensional
right $A$-module. Let $m$  be a basis of the module
$\Ext_{A}^{p+s}(_A\kk,{_AA})$.  It follows from  the right version
of \cite[Lemma 2.13 (1)]{rs} that $ m\cdot x_i=0$ for all $1\se i\se
\tt$. Since $\Ext_{C}^{p+s}(_C\kk,{_CC})\cong\Ext_{A}^{p+s}(_A\kk,{_AA})$ as
right $\bgm$-modules, we have showed that $$\Ext_{A}^i(_A\kk,{_AA})\cong\begin{cases}0,&i\neq p+s;\\
\kk_\xi ,&i=p+s.\end{cases}$$

Similarly, we have
$$\dim\Ext^i_A(\kk_A,A_A)=\begin{cases}0,& i\neq p+s;
\\1,&i=p+s.\end{cases}$$
By Lemma \ref{global}, the algebra $A$ is AS-regular of global
dimension $p+s$. \qed

Now we can give  a necessary and sufficient condition for a pointed
Hopf algebra $U(\mc{D},\lmd)$ to be CY.

\begin{thm}\label{cy}
Let $(\mc{D},\lmd)$ be a generic datum of finite Cartan type  and $A$
the Hopf algebra $U(\mc{D},\lmd)$. Let $s$ be the rank of $\bgm$ and $p$
the number of the positive roots of the Cartan matrix.
\begin{enumerate}
\item[(a)] The rigid dualizing complex of the Hopf algebra
$A=U(\mc{D},\lmd)$ is $_\psi A[p+s]$, where $\psi$ is defined by $
\psi(x_k)=\prod_{i=1,i\neq j_k}^{p}\chi_{_{\bt_i}}(g_k)x_k$, for all
$1\se k\se \tt$, and $\psi(g)=(\prod_{i=1}^p\chi_{_{\bt_i}})(g)$ for
all $g\in \bgm$, where each $ j_k $, $1\se k\se \tt$,  is the integer
such that $\bt_{j_k}=\al_k$.
\item[(b)]  The algebra $A$ is CY if and only if
$\prod_{i=1}^p\chi_{_{\bt_i}}=\varepsilon$ and $\mc{S}^2_A$ is an
inner automorphism.
\end{enumerate}
\end{thm}
\proof (a) By \cite[Prop. 4.5]{bz} and Theorem \ref{as}, the
rigid dualizing complex of $A$ is isomorphic to
$_{[\xi]\mc{S}^2_A}A[p+s]$, where $\xi$ is the algebra homomorphism
defined in Theorem \ref{as}. It is not difficult to see that
$$([\xi]\mc{S}^2_A)(g)=(\prod_{i=1}^p\chi_{_{\bt_i}})(g),$$ for all $g\in \bgm$. For $1\se k\se \tt$, we have $ \Delta (x_k)=x_k\ot 1+g_k\ot
x_k$ and $\mc{S}^2_A(x_k)=\chi_k(g_k^{-1})x_k$. If $j_k$ is the
integer such that $\bt_{j_k}=\al_k$, then
$\chi_{_{\bt_{j_k}}}(g_k)=\chi_k(g_k)$. So
$$\begin{array}{ccl}([\xi]\mc{S}^2_A)(x_k)&=&\chi_k(g_k^{-1})[\xi](x_k)\\&=&\chi_k(g_k^{-1})\prod_{i=1}^{p}\chi_{_{\bt_i}}(g_k)(x_k)\\&=&\prod_{i=1,i\neq
j_k}^{p}\chi_{_{\bt_i}}(g_k)(x_k).\end{array}$$

(b) follows from Theorem \ref{as}  and \cite[Thm. 2.3]{hvz}. \qed

\begin{rem}\label{lifting} From Theorem \ref{cy}, we can see that for a pointed Hopf algebra $U(\mc{D},\lmd)$, it is CY \iff\  its associated graded algebra $U(\mc{D},0)$ is
CY.
\end{rem}



\begin{cor}
Assume that $A=U(\mc{D},\lmd)$. For every $A$-$A$-bimodule $M$,
there are isomorphisms:
\begin{equation}\label{vdbpointed}\HH^i(A,M)\cong
\HH _{p+s-i}(A,{_{\psi^{-1}}M}),\;\;0\se i\se p+s,\end{equation}
where $\psi$ is the algebra automorphism defined in Theorem
\ref{cy}.
\end{cor}
\proof This follows from \cite[Cor. 5.2]{bz} and Theorem
\ref{as}.\qed

\newpage
\section{Calabi-Yau Nichols algebras of finite Cartan type}

As we noted in Remark \ref{lifting}, the CY property of $U(\mc{D},\lmd)$ is determined by the CY property of $U(\mc{D},0)$, which is equal to the smash product $\mc{B}(V)\# \kk\bgm$ of the Nichols algebra $\mc{B}(V)$ with the group algebra $\kk\bgm$. One would naturally ask whether or not the CY property of $U(\mc{D}, 0)$ depends on the CY property of $\mc{B}(V)$. In this section, we work out a criterion for the Nichols algebra $\mc{B}(V)$ to be Calabi-Yau, and answer the question in Section 4. We fix a  generic datum of finite Cartan type
$$(\mc{D},0)=(\bgm,(a_{ij})_{1\se i,j\se \tt},(q_{_I})_{I\in \mathcal {X}},
(g_i)_{1\se i\se \tt},(\chi_i)_{1\se i\se \tt},0),$$ where $\bgm$  is a free abelian group of
rank $s$.  Let $V$ be the generic
braided vector space with basis $\{x_1,\cdots,x_\tt\}$ whose
braiding is given by
$$c(x_i\ot x_j)=q_{ij}x_j\ot x_i$$ for  all $1\se i,j\se \tt$, where
$q_{ij}=\chi_j(g_i)$. Recall that the Nichols algebra $\mathcal
{B}(V)$ is generated by $x_i$, $1\se i\se \tt$, subject to the
relations
$$\ad_c(x_i)^{1-a_{ij}}x_j= 0,\;\; 1\se i,j\se \tt,\;\;i\neq j,$$
where $\ad_c$ is the braided adjoint representation.

By \cite[Thm. 4.3]{as3}, the Nichols algebra $\mc{B}(V)$ is a
subalgebra of $U(\mc{D},0)$, and the following monomials in  root
vectors
$$\{x_{_{\bt_1}}^{a_1}\cdots x_{_{\bt_p}}^{a_p}\mid a_i\le 0,1\se i\se p\}$$ form a PBW basis of the Nichols algebra $\mc{B}(V)$.  The degree
(cf. (\ref{deg pointd})) of each PBW basis
element  is defined by $$\deg (x_{_{\bt_1}}^{a_1}\cdots
x_{_{\bt_p}}^{a_p})=(a_1,\cdots, a_p,\sum a_iht(\bt_i))\in
(\mathbb{Z}^{\le0})^{p+1},$$ where $ht(\bt_i)$ is the height of
$\bt_i$.

The following result is a direct consequence of Lemma \ref{filt1}.

\begin{lem}\label{nichols filt}
In the Nichols algebra $\mc{B}(V)$, for $j>i$, we have
\begin{equation}\label{comm1}[x_{_{\bt_i}}, x_{_{\bt_j}}]_c=\sum_{\textbf{a}\in \mathbb{N}^p}\rho_{\textbf{a}}x_{_{\bt_1}}^{a_1}\cdots x_{_{\bt_p}}^{a_p}, \end{equation}
where $\rho_{\textbf{a}}\in \kk$ and $\rho_{\textbf{a}}\neq0$ only
when $\textbf{a}=(a_1,\cdots,a_p)$ satisfies that $a_k=0$ for $k\se
i$ and $k\le j$.
\end{lem}

Order the PBW basis elements by degree as in (\ref{ordering}).  By
Lemma \ref{nichols filt}, we obtain the following corollary.

\begin{cor}\label{nichols gra}
The Nichols algebra $\mc{B}(V)$ is an $\NN^{p+1}$-filtered algebra,
whose associated graded algebra $\mb{G}\rr\mc{B}(V)$ is isomorphic
to the following algebra:
$$\kk\lan x_{_{\bt_1}},\cdots,x_{_{\bt_p}}\mid x_{_{\bt_i}}x_{_{\bt_j}}=\chi_{_{\bt_j}}(g_{_{\bt_i}})x_{_{\bt_j}}x_{_{\bt_i}},\;\;1\se
i<j\se p\ran,$$ where $x_{_{\bt_1}},\cdots,x_{_{\bt_p}}$ are the
root vectors of $\mc{B}(V)$.
\end{cor}

For elements $\{x_{_{\bt_1}}^{a_1}\cdots x_{_{\bt_p}}^{a_p}\},$
where $a_1,\cdots, a_p\le 0$, define \vspace{3mm}\\\centerline{$
d_0(x_{_{\bt_1}}^{a_1}\cdots x_{_{\bt_p}}^{a_p})=\sum_{i=1}^p
a_iht(\bt_i).$}\vspace{3mm}\\Then $R=\mc{B}(V)$ is a graded algebra
with grading given by $d_0$. Let $R^{(0)}=R$. Define
$d_1(x_{_{\bt_1}}^{a_1}\cdots x_{_{\bt_p}}^{a_p})=a_p$. We obtain an
$\NN$-filtration on $R^{(0)}$. Let $R^{(1)}=\Gr R^{(0)}$ be the
associated graded algebra. In a similar way, we define
$d_2(x_{_{\bt_1}}^{a_1}\cdots x_{_{\bt_p}}^{a_p})=a_{p-1}$ and let
$R^{(2)}=\Gr R^{(1)}$ be the associated graded algebra. Inductively,
we obtain a sequence of $\NN$-filtered algebras $R^{(0)},\cdots,
R^{(p)}$, such that $R^{(i)}=\Gr R^{(i-1)}$, for $1\se i\se p$, and
$R^{(p)}=\mb{G}\rr R$.


The algebra $R^e$ has a PBW basis as follows
\begin{equation}\label{pbw2}\{x_{_{\bt_1}}^{a_1}\cdots
x_{_{\bt_p}}^{a_p}\ot x_{_{\bt_p}}^{b_p}\star\cdots \star
x_{_{\bt_1}}^{b_1}\mid a_1,\cdots, a_p,b_1, \cdots, b_p\le 0
\},\end{equation}where  ``$\star$'' denotes the multiplication in
$R^{op}$. Similarly, define a degree on each element   as
$$\begin{array}{cl}&\deg (x_{_{\bt_1}}^{a_1}\cdots x_{_{\bt_p}}^{a_p}\ot
x_{_{\bt_p}}^{b_p} \star\cdots \star
x_{_{\bt_1}}^{b_1})\\=&(a_1+b_1, \cdots, a_p+b_p,  \sum
(a_i+b_i)ht\bt_i)\in (\ZZ^{\le0})^{(p+1)}.\end{array}$$ Then $R^e$
is an $\NN^{p+1}$-filtered algebra whose associated graded algebra
$\mb{G}\rr(R^e)$ is isomorphic to $(\mb{G}\rr R)^e$.

In a similar way, we obtain a sequence of $\NN$-filtered
algebras $(R^e)^{(0)},\cdots, (R^e)^{(p)}$, such that
$(R^e)^{(i)}=\Gr((R^e)^{(i-1)})$, for $1\se i\se p$, and
$(R^e)^{(p)}=\mb{G}\rr R^e$. In fact, $(R^e)^{(i)}=(R^{(i)})^e$, for
$0\se i\se p$.

\begin{lem}\label{noetherian}
Let $R=\mc{B}(V)$ be the Nichols algebra of $V$. Then the algebra
$R^e$ is Noetherian.
\end{lem}
\proof The sequence $(R^e)^{(0)},\cdots, (R^e)^{(p)}$ is a sequence
of algebras, each of which is the associated graded algebra of the
previous one with respect to an $\NN$-filtration. The algebra
$(R^e)^{(p)}$ is isomorphic to $(\mb{G}\rr R)^e$, which is
Noetherian. By \cite[Thm. 1.6.9]{mr}, the algebra $R^e$ is
Noetherian. \qed

\begin{lem}\label{lem hs}
The algebra $R=\mc{B}(V)$ is homologically smooth.
\end{lem}
\proof  Since $R^e$ is Noetherian by Lemma \ref{noetherian} and $R$
is a finitely generated $R^e$-module, it is sufficient to  prove
that the projective dimension $\projdim{}_{R^e}R$ is finite.  The
filtration on each $(R^{(i)})^e$, $0\se i\se p-1$, is bounded below.
In addition, from the proof of the foregoing Lemma \ref{noetherian},
each $(R^{(i)})^e$ is Noetherian for $0\se i\se p$. Therefore,
$(R^{(i)})^e$ is a Zariskian algebra for each $0\se i\se p-1$. It is
clear that each $R^{(i)}$, $1\se i\se p-1$, viewed as an
$(R^{(i)})^e$-module has a good filtration. By \cite[Cor. 5.8]{lv1},
we have
$$\begin{array}{ccl}\projdim_{R^e}R&=&\projdim_{(R^{(0)})^e}R^{(0)}\se\projdim_{(R^{(1)})^e}R^{(1)}\se\cdots\\
&\se & \projdim_{(R^{(p)})^e}R^{(p)}=\projdim_{(\mb{G}\rr
R)^e}\mb{G}\rr R.\end{array}$$ The algebra $\mb{G}\rr R$ is a
quantum polynomial algebra of $p$ variables.  From the Koszul
bimodule complex of $\mb{G}\rr R$ (cf. (\ref{kosbimod})), we obtain
that $\projdim_{(\mb{G}\rr R)^e}\mb{G}\rr R=p$. Therefore,
$\projdim_{R^e}R\se p$ and $R$ is homologically smooth. \qed

\begin{prop}\label{nicholsas}
Let $R=\mc{B}(V)$ be the Nichols algebra of $V$.
\begin{enumerate}\item[(a)] $R$ is AS-regular of global dimension
$p$.
\item[(b)]  The rigid dualizing complex of $R$ in the  graded sense  is
isomorphic to  $ _\vph R(l)[p]$ for some integer  $l$ and some
$\NN$-graded algebra automorphism $\vph$ of degree 0.
\item[(c)] The rigid dualizing complex in the ungraded sense is just
$ _\vph R[p]$.
\end{enumerate}
\end{prop}
\proof Let $x_{_{\bt_1}},\cdots,x_{_{\bt_p}}$ be the root vectors.
By Lemma \ref{nichols filt}, we can use a similar argument to the proof of Lemma \ref{global} to show that the algebra $R$ is an iterated graded Ore extension of $\kk[x_{_{\bt_1}}]$. Indeed,
$$R\cong \kk[x_{_{\bt_1}}][x_{_{\bt_2}};\tau_2,\delta_2]\cdots[x_{_{\bt_p}};\tau_p,\delta_p],$$
where for $2\se j\se p$, $\tau_j$ is an algebra automorphism such
that $\tau_{j}(x_{_{\bt_i}})$ is just a scalar multiple of
$x_{_{\bt_i}}$ for $i<j$, and $\delta_j$ is a $\tau_{j}$-derivation
such that $\delta_j(x_{_{\bt_i}})$, $i<j$, is a linear combination
of monomials in $x_{_{\bt_{i+1}}},\cdots, x_{_{\bt_{j-1}}}$.
 It is
well-known that $\kk[x_{_{\bt_1}}]$ is an AS-regular algebra of
dimension 1, and the AS-regularity is preserved under graded Ore
extension. So $R$ is an AS-regular algebra of dimension $p$.
Therefore,  the rigid dualizing complex of $R$ in the graded case is
isomorphic to $ _\vph R(l) [p]$ for some graded algebra automorphism
$\vph$ and some $l\in \mathbb{Z}$. By Lemma \ref{noetherian}, $R^e$
is Noetherian. Thus the rigid dualizing complex $_\vph R(l)[p]$ in
the graded case implies the rigid dualizing complex $_\vph R[p]$ in the
ungraded case. \qed

We claim that the automorphism $\vph$ in Proposition \ref{nicholsas}
is just a scalar multiplication. We need some preparations to prove
this claim.

If $R$ is a $\bgm$-module algebra, then the algebra $R^e$ is also a
$\bgm$-module algebra with the natural action $g(r\ot s):={g(r)}\ot
{g(s)}$, for all $g\in \bgm$ and $r,s\in R$.

\begin{lem}\label{gaction}
Let $R$ be a $\bgm$-module algebra, such that $\kk^\times$ is the
group of units of $R$. Assume that $U$ is an $R^e\#\kk\bgm$-module,
and  $U\cong { R_\phi}$, for an algebra automorphism $\phi$, as $R^e\# \kk\bgm$-modules. Then
\begin{enumerate}
\item[(a)] the algebra automorphism $\phi$
preserves $\bgm$-action;
\item[(b)] the  $R^e\#\kk\bgm$-module
structures on $U$ (up to isomorphism) are parameterized  by
$\Hom(\bgm,\kk)$, the set of group homomorphisms from $\bgm$ to
$\kk^\times$.
\end{enumerate}
\end{lem}
\proof  (a) Fix an isomorphism $U\cong { R_\phi}$. Let $u\in U$ be the
element mapped to $1\in R$. Then $U=Ru$ and we have $g(ru)={
g(r)}g(u)$ for all $r\in R$ and $g\in \bgm$. To determine the
$\bgm$-action on $U$, we only need to determine $g(u)$ for $g\in
\bgm$. Since $g(u)\in U$,  there is some $r_g\in R$, such that
$g(u)=r_gu$. On the other hand, we have
$$U=g(U)=g(Ru).$$ So there is some $s\in R$, such that
$u={g(s)}r_gu$. Since the element $u$ forms an $R$-basis of $U$, the
element $r_g$ has a left inverse. Similarly, there is some $s'\in
R$, such that $u=r_gug(s')$.  Since $U\cong { R_\phi}$ as
$R$-$R$-bimodules, we have
\begin{equation}\label{equa 1}\phi(r)u=ur,\end{equation} for any $r\in R$. So
$u=r_gug(s')=r_g\phi^{-1}(g(s'))u$. Thus $r_g$ has a right inverse
as well.  Consequently, $r_g$ is a unit in $R$, and $r_g\in
\kk^\times$. We have $g(h(u))=(gh)(u)$ for $g,h\in \bgm$. That
is, $r_{gh}=r_gr_h$. Therefore, the $\bgm$-action on $U$ defines a
group homomorphism from $\bgm$ to $\kk^\times$, denoted by
$\chi:\bgm\ra \kk^\times$.  Since $U$ is an
$R^e\#\kk\bgm$-module, we have $g(rus)={g(r)}{g(u)}{g(s)}$,
for any $r,s\in R$ and $g\in \bgm$. To show that $\phi$ preserves the $\bgm$-action, we compute
$g(\phi(r)(u))$. On one hand,
 we have
$$\begin{array}{ccl}{g(\phi(r)u)}&=&{g(ur)}\\&=&g(u)g(r)\\&=&\chi(g)u{g(r)}\\&\overset{\tiny(\ref{equa 1})}=&\chi(g)\phi(g(r))u.\end{array}$$
On the other hand, we have
$$\begin{array}{ccl}{g(\phi(r)u)}&=&{g(\phi(r))}g(u)\\&=&\chi(g){g(\phi(r))}u.\end{array}$$ So ${g(\phi(r))}=\phi({g(r)})$. That is, the automorphism $\phi$ preserves $\bgm$-action.

(b) In Part (a) we have shown that a $\bgm$-action on $U$ is determined by a group homomorphism from $\bgm$ to $k^\times$ such
that $U$ is an $R^e\#\kk\bgm$-module.

Suppose there are two $\bgm$-actions on $U$ such that they are
isomorphic. We write these two actions as ${g^{\cdot1}}(u)=r_gu$ and
${g^{\cdot2}}(u)=s_gu$. Denote by  $U_1$ and $U_2$  the
$\bgm$-modules with these two actions respectively. Let $f:U_1\ra
U_2$ be an $R^e\#\kk\bgm$-module isomorphism. Then $f(u)=ru$ for
some unit $r\in R$. Since the set of units of $R$ is $\kk^\times$,
we have $r\in \kk^\times$. On one hand, we have
$$\begin{array}{ccl}f({g^{\cdot1}}(u))&=&f(r_gu)\\&=&r_gr u.\end{array}$$
On the other hand, we also have
$$\begin{array}{ccl}f({g^{\cdot1}}(u))&=&{g^{\cdot2}}(f(u))\\&=&{g^{\cdot2}}(ru)\\&=&r{g^{\cdot2}}(u)\\&=&s_gru.
\end{array}$$
Therefore, $r_g=s_g$, and (b) follows. \qed

If $U$ is an $R^e\#\kk\bgm$-module, then we can define an
$(R\#\kk\bgm)^e$-module $U\#\kk\bgm$. It is isomorphic to $U\ot \bgm$
as vector space with bimodule structure given by
$$(r\#h)(u\ot g):=rh(u)\ot hg,$$
$$(u\ot g)(r\#h):=u{g(r)}\ot gh,$$ for any $r\#h\in R\#H$ and $u\ot g\in U\ot \bgm$.

\begin{lem}\label{nichols bimod}
Let $R$ be a $\bgm$-module algebra with $\kk^\times$ being the group
of units and $U$ an $R^e\#\kk\bgm$-module. Assume that  $U\cong {
R_\phi}$ as $R^e\#\kk\bgm$-modules, where $\phi$ is an algebra
automorphism. If the $\bgm$-action on $U$ is defined by a group
homomorphism $\chi:\bgm\ra \kk^\times$. Then $U\#\kk\bgm\cong
{(R\#\kk\bgm)_\psi}$ as $(R\#\kk\bgm)^e$-modules, where $\psi$ is
the algebra automorphism defined by
$\psi(r\#g)=\chi(g^{-1})\phi(r)\#g$ for any $r\#g\in R\#\kk\bgm$.
\end{lem}
\proof The homomorphism $\psi$ defined in the lemma is clearly
bijective. First we check that it is an  algebra homomorphism. For any
$r\#g,s\#h\in R\#\kk\bgm$, we have
$$\begin{array}{ccl}\psi((r\#g)(s\#h))&=&\psi(rg(s)\#gh)\\&=&\chi(h^{-1}g^{-1})\phi(rg(s))\#gh\\
&=&\chi(h^{-1}g^{-1})\phi(r)\phi(g(s))\#gh\\
&=&\chi(h^{-1}g^{-1})\phi(r)g(\phi(s))\#gh\\
&=&(\phi(r)\chi(g^{-1})\#g)({\phi(s)}\chi(h^{-1})\#h)\\
&=&\psi(r\#g)\psi(s\#h).
\end{array}$$ The forth equation holds since $\phi$ preserves the $\bgm$-action by Lemma
\ref{gaction}.

Next we show that $U\#\kk\bgm\cong {(R\#\kk\bgm)_\psi}$ as
$(R\#\kk\bgm)^e$-modules. Fix an isomorphism $U\cong { R_\phi}$ and
let $u\in U$ be
 the element mapped to $1\in R$. We define a
homomorphism $\Phi:U\#\kk\bgm\ra {(R\#\kk\bgm)_\psi}$ by $\Phi(ru\ot
g)=\chi(g^{-1})r\#g$. It is easy to see that $\Phi$ is an
isomorphism of left $R\#\kk\bgm$-modules. Now we show that it is a
right $R\#\kk\bgm$-module homomorphism. Indeed, we have
$$\begin{array}{ccl}\Phi(u(r\#g))&=&\Phi(ur\ot g)\\
&=&\Phi(\phi(r)u\ot g)\\
&=& \chi(g^{-1})\phi(r)\#g\\
&=& \Phi(u)\psi(r\#g) \\
&=&   \Phi(u)\cdot(r\#g)  .\end{array}$$ \qed

Now we can prove the following lemma.

\begin{lem}\label{scalar}
Keep the notations as in  Proposition \ref{nicholsas}. The actions
of $\vph$ on generators $x_1,\cdots, x_\tt$ are just scalar
multiplications. \end{lem}

\proof By Proposition \ref{nicholsas} and Lemma \ref{vdb}, we have
$R$-$R$-bimodule isomorphisms
$$\Ext_{R^e}^i(R,R^e)\cong \begin{cases}0,& i\neq p;
\\ R_\vph   ,&i=p.\end{cases}$$

The group $\bgm$ is a free abelian group of rank $s$, so the algebra
$\kk\bgm$ is a CY algebra of dimension $s$. Following from \cite[Sec. 2]{f}, $R_\vph$
is an $R^e\#\kk\bgm$-module and there are $(R\#\kk\bgm)^e$-bimodule
isomorphisms
$$\Ext_{(R\#\kk\bgm)^e}^i(R\#\kk\bgm,(R\#\kk\bgm)^e)\cong\begin{cases}0,& i\neq p+s;
\\( R_\vph)\#\kk\bgm,&i=p+s.\end{cases}$$
For the sake of completeness, we sketch the proof here.
  By Lemma \ref{lem hs}, $R$ is homologically
smooth. That is, $R$ has a bimodule projective resolution
\begin{equation}\label{equ 4}0\ra P_q\ra \cdots \ra P_1\ra P_0\ra R\ra 0,\end{equation} with each $P_i$
being finitely generated as an $R$-$R$-bimodule.

$\Ext_{R^e}^*(R,R^e)$ are the cohomologies of the complex
$\Hom_{R^e}(P_\bullet,R^e)$.
 The algebra $R^e$ is an $R^e\#\kk\bgm$-module defined by $$((c\ot
d)\#g)\cdot (a\ot b)=g(a)d\ot cg(b)$$ for any $a\ot b\in R^e$ and
$(c\ot d)\#g\in R^e\#\kk\bgm$. Then each $\Hom_{R^e}(P_i,R^e)$ is an
 $R^e\#\kk\bgm$-module as well:
\begin{equation}\label{equ 3}[((c\ot d)\#g)\cdot f](x)=((c\ot d)\#g)\cdot f(x),\end{equation}where $(c\ot d)\#g\in
R^e\#\kk\bgm$, $f\in \Hom_{R^e}(P_i,R^e)$ and $x\in P_i$. Now $\Hom_{R^e}(P_\bullet,R^e)$ is a complex
of left $R^e\#\kk\bgm$-modules. Thus we obtain that
$\Ext_{R^e}^p(R,R^e)\cong R_\vph$ is an $R^e\#\kk\bgm$-module.

Let $A=R\#\kk\bgm$. Observe that $A^e$ is an
$R^e\#\kk\bgm$-$A^e$-bimodule. The left $\kk\bgm$-module action is
defined by \begin{equation}\label{equ 2}g\cdot(a\#h\ot
b\#k)=g(a)gh\ot b\#kg^{-1},\end{equation} for any $a\#h\ot b\#k\in
A^e$ and $g\in \bgm$. The left $R^e$-action and right $A^e$-action
are given by multiplication. Let $W$ be the vector space $\kk\bgm\ot
\kk\bgm$. Then  $R^e\ot W$ is also an $R^e\#\kk\bgm$-$A^e$-bimodule
defined by
$$((c\ot d)\#g)\cdot(a\ot b\ot h\ot k)=cg(a)\ot g(b)d\ot gh\ot kg^{-1}$$
and $$(a\ot b\ot h\ot k)\cdot(c\#h'\ot d\#k')=ah(c)\ot
((k^{-1}k'^{-1})d)b \ot hh'\ot k'k.$$ It is not difficult to see
that the morphism $f:A^e\ra R^e\ot W$ defined by
$$f(a\#h\ot b\#k)=a\ot k^{-1}(b) \ot h\ot k$$ is an isomorphism of
$R^e\#\kk\bgm$-$A^e$-bimodules.

Let $P$ be a finitely generated projective $R^e$-module. The
$\kk\bgm$-$A^e$-bimodule structure of $R^e\ot W$  induces a
$\kk\bgm$-$A^e$-bimodule structure on $\Hom_{R^e}(P,R^e\ot W)$. We
define a $\kk\bgm$-$A^e$-bimodule structure on $\Hom_{R^e}(P,R^e)\ot W$ as
follows
$$g\cdot (f\ot h\ot k)=g\cdot f\ot gh\ot kg^{-1}$$ and
$$(f\ot h\ot k)\cdot(c\#h'\ot d\#k')= (h(c)\ot (k^{-1}k'^{-1})d)\cdot f\ot  hh'\ot
k'k,$$ where the $R^e\#\kk\bgm$-module structure on
$\Hom_{R^e}(P,R^e)$ is defined in (\ref{equ 3}). Now the canonical
isomorphism from $\Hom_{R^e}(P,R^e)\ot W$ to $\Hom_{R^e}(P,R^e\ot
W)$ is a $\kk\bgm$-$A^e$-bimodule isomorphism.

Since  $R$ admits a resolution like
(\ref{equ 4}) with each $P_i$ finitely generated,  we have
$$\Ext^i_{R^e}(R,R^e\ot W)\cong \Ext^i_{R^e}(R,R^e)\ot W$$ as
$\kk\bgm$-$A^e$-bimodules for all $i\le 0$. On the other hand,
we have Stefan's spectral sequence \cite{st}:
$$\Ext_{\kk\bgm}^m(\kk,\Ext_{R^e}^n(R,A^e))\Rightarrow \Ext_{A^e}^{m+n}(A,A^e).$$ 
Thus for $m,n\le0$, we have $$\begin{array}{ccl}\Ext_{\kk\bgm}^{m}(\kk,\Ext^n_{R^e}(R,A^e))&\cong&\Ext_{\kk\bgm}^{m}(\kk,\Ext^n_{R^e}(R,R^e\ot W))\\
&\cong&\Ext_{\kk\bgm}^{m}(\kk,\Ext^n_{R^e}(R,R^e)\ot W).
\end{array}$$
Hence,  $\Ext_{\kk\bgm}^m(\kk,\Ext_{R^e}^n(R,A^e))=0$ except that $m=s$
and $n=p$. Therefore,
$$\Ext_{(R\#\kk\bgm)^e}^i(R\#\kk\bgm,(R\#\kk\bgm)^e)=0$$ for $i\neq
p+s$ and
$$\Ext_{A^e}^{p+s}(A,A^e)\cong\Ext_{\kk\bgm}^{s}(\kk,\Ext^p_{R^e}(R,A^e)).$$

Let $M$ be a left $\kk\bgm$-module. One can consider it as a
$\kk\bgm$-$\kk\bgm$-bimodule $M_\vps$ with the trivial right
$\kk\bgm$-module action. The algebra $\kk\bgm$ is a CY algebra of
dimension $s$. From Van den Bergh's duality theorem (Theorem
\cite[Thm. 1]{v1}) we obtain the following isomorphisms:
\begin{equation}\label{equ 1}\begin{array}{ccl}\Ext^s_{\kk\bgm }(\kk,M)&\cong&\HH^s(\kk\bgm,M_\vps)\\&\cong&\HH_0(\kk\bgm,M_\vps)\\&\cong&\Tor^{\kk\bgm}_0(\kk,M).\end{array}\end{equation}
Now we have the following isomorphisms of right $A^e$-modules
$$\begin{array}{ccl}\Ext_{A^e}^{p+s}(A,A^e)&\cong&\Ext_{\kk\bgm}^{s}(\kk,\Ext^p_{R^e}(R,A^e))\\
&\cong&\Ext_{\kk\bgm}^{s}(\kk,\Ext^p_{R^e}(R,R^e)\ot W)\\
&\cong&\Ext_{\kk\bgm}^{s}(\kk,R_\vph\ot W)\\
&\cong&\Tor^{\kk\bgm}_0(\kk,R_\vph\ot W)\\
&\cong& \kk\ot _{\kk\bgm}R_\vph\ot W.\end{array}$$
If we look at the $\kk\bgm$-$A^e$-bimodule structure  on $R_\vph\ot W$ carefully, we obtain that $$\kk\ot _{\kk\bgm}R_\vph\ot W\cong R_\vph\#\kk\bgm$$ as right $A^e$-modules.

Since the connected graded algebra $R$ is a domain by Lemma
\ref{pointeddef}, the group of the units of $R$ is $\kk^\times$.
Following Lemma \ref{gaction} and \ref{nichols bimod}, we have $(
R_\vph)\#\kk\bgm\cong {(R\#\kk\bgm)_{\bar{\psi}}}$, where
$\bar{\psi}$ is the algebra automorphism  defined by
$\bar{\psi}(r\#g)={\vph}(r)\chi(g^{-1})$ for some algebra
homomorphism $\chi:\bgm\ra \kk$.

On the other hand,  since  $A=R\#\kk\bgm\cong U(\mc{D},0)$,  we have
$A$-$A$-bimodule isomorphisms:
$$\Ext_{A^e}^i(A,A^e)\cong\begin{cases}0,& i\neq p+s;
\\A_\psi,&i=p+s,\end{cases}$$
where $\psi$ is the algebra automorphism defined in Theorem
\ref{cy}.

Therefore, we obtain an $A$-$A$-bimodule isomorphism
$A_{\bar{\psi}}\cong A_\psi.$ That is, $\bar{\psi}$ and $\psi$
differ only by an inner automorphism. By Lemma \ref{pointeddef},
the graded algebra $A$ is a domain; and the invertible elements of $A$
are in $\kk \bgm$. The actions of $\psi$ and the group actions on
generators $x_1,\cdots, x_\tt$ are just scalar multiplications. Thus
the actions of $\bar{\psi}$ on $x_1,\cdots, x_\tt$ are also scalar
multiplications. Since $\bar{\psi}(x_i)=\vph(x_i)$ for all $1\se
i\se \tt$,  we obtain the desired result.\qed

Now we are ready to  prove the main theorem of this section.

\begin{thm}\label{nichols cy}
Let $V$ be a generic braided vector space of finite Cartan type, and
$R=\mc{B}(V)$ the Nichols algebra of $V$. Let $p$ be the number of the positive roots of the Cartan matrix. For each $1\se k\se \tt$,
let $ j_k $ be the integer  such that $\bt_{j_k}=\al_k$.
\begin{enumerate}
\item[(a)] The rigid dualizing complex is isomorphic to $  _\vph R[p]$, where $\vph$ is the  algebra automorphism  defined by $$\vph(x_k)=(\prod_{i=1}^{j_k-1}\chi^{-1}_k(g_{_{\bt_i}}))(
\prod_{i=j_k+1}^{p}\chi_{_{\bt_i}}(g_k))x_k,$$ for all $1\se k\se
\tt$.
\item[(b)] The algebra $R$ is a CY algebra if and only if $$\prod_{i=1}^{j_k-1}\chi_k(g_{_{\bt_i}})= \prod_{i=j_k+1}^{p}\chi_{_{\bt_i}}(g_k),$$ for all $1\se k\se \tt$.
\end{enumerate}
\end{thm}
\proof

(a)  Note that $\mb{G}\rr R$ is isomorphic to the following quantum
polynomial algebra:
$$\kk\lan x_{_{\bt_1}},\cdots,x_{_{_{\bt_p}}}\mid x_{_{\bt_i}}x_{_{\bt_j}}=\chi_{_{\bt_j}}(g_{_{\bt_i}})x_{_{\bt_j}}x_{_{\bt_i}},\;\;1\se
i<j\se p\ran.$$ By \cite[Prop. 8.2 and Thm. 9.2]{vdb}, $\mb{G}\rr R$ has a rigid
dualizing complex $_{\bar{\zeta}}\mb{G}\rr R[p](\cong \mb{G}\rr
R_{\bar{\zeta}^{-1}}[p]$), where $\bar{\zeta}$ is defined by
$$\bar{\zeta}(x_{_{\bt_k}})=\chi^{-1}_{_{\bt_k}}(g_{_{\bt_1}})\cdots\chi^{-1}_{
{k}}(g_{_{\bt_{k-1}}})\chi_{_{\bt_{k+1}}}(g_{_{\bt_k}})\cdots\chi_{_{\bt_p}}(g_{_{\bt_k}})x_{_{\bt_k}},$$
for all $1\se k\se p$.

On the other hand, it follows from Proposition \ref{nicholsas} and
Lemma \ref{scalar} that $R$ has a rigid dualizing complex
$_{\vph}R$, where
 $\vph$ is an algebra automorphism such that for each  $1\se
k\se \tt$,  $\vph(x_k)$ is a scalar multiple of $x_k$.  Assume that
$\vph(x_k)=l_kx_k$, with $l_k\in \kk$.

Let $R^{(0)},\cdots,R^{(p)}$ be the sequence of algebras defined
after Corollary \ref{nichols gra}. By Lemma \ref{nichols filt},
and applying a similar argument to the one  in the proof of Proposition
\ref{nicholsas}, we obtain that each $R^{(i)}$, $0\se i\se p$, is an
iterated Ore extension of the polynomial algebra $\kk[x]$. Thus each
of them is AS-regular. It follows from \cite[Prop. 1.1]{y2} that
each $R^{(i)}$, $1\se i\se p$, has a rigid dualizing complex $
_{\vph^{(i)}}(R^{(i)})[p]$, where ${\vph^{(i)}}=\Gr{\vph^{(i-1)}}$
and $\vph^{(0)}=\vph$. Since for each $1\se k\se \tt$,
$\vph(x_k)=l_kx_k$, we have ${\vph^{(p)}}(x_k)=l_kx_k$. Because
$R^{(p)}=\mb{G}\rr R$, there is a bimodule isomorphism  $
_{\vph^{(p)}}(R^{(p)})\cong {}_{\bar{\zeta}}(\mb{G}\rr R).$ We
obtain that ${\vph^{(p)}} =\bar{\zeta}$, as $R$ is connected.
Therefore, for each $1\se k\se \tt$,
$$l_kx_k=\bar{\zeta}(x_k)=(\prod_{i=1}^{j_k-1}\chi^{-1}_k(g_{_{\bt_i}}))(
\prod_{i=j_k+1}^{p}\chi_{_{\bt_i}}(g_k))x_k,$$ where $j_k$ is the
integer such that $\bt_{j_k}=\al_k$.

Now we can conclude that
$\vph(x_k)=(\prod_{i=1}^{j_k-1}\chi^{-1}_k(g_{_{\bt_i}}))(
\prod_{i=j_k+1}^{p}\chi_{_{\bt_i}}(g_k))x_k,$ for each $1\se k\se
\tt$.

(b) The algebra $R$ is homologically smooth by Lemma \ref{lem hs}.
It follows from Corollary \ref{cor cyrid} that $R$ is CY if and only if $R\cong
{}_\vph R$ as bimodules. That is, $R$ is CY if and only if
$\vph=\id$. Hence (b) follows from (a).
 \qed

\begin{eg}\label{eg 1}
Let $(\mc{D},\lmd)=(\bgm,(a_{ij}),(q_{_I}),(g_i),(\chi_i),0)$ be a generic datum such
that the Cartan matrix is of type  $A_2$. This defines a braided
vecter space $V$. Let $\{x_1,x_2\}$ be a basis of $V$. The braiding
of $V$ is given by
$$c(x_i\ot x_j)=\chi_j(g_i)x_j\ot x_i,\;\;i,j=1,2.$$
The Nichols algebra $R=\mc{B}(V)$ of $V$ is generated by $x_1$ and
$x_2$ subject to the relations
$$x_1^2x_2-q_{_{12}}x_1x_2x_1-q_{_{11}}q_{_{12}}x_1x_2x_1+q_{_{11}}q_{_{12}}^2x_2x_1^2=0,$$
$$x_2^2x_1-q_{_{21}}x_2x_1x_2-q_{_{22}}q_{_{21}}x_2x_1x_2+q_{_{22}}q_{_{21}}^2x_1x_2^2=0,$$
where $q_{ij}=\chi_j(g_i)$.
 The element $s_1s_2s_1$ is the longest element in the Weyl group $\mc{W}$. Let $\al_1$ and $\al_2$ be the two simple roots. Then the positive roots are as follows $$\bt_1=\al_1,\;\;\bt_2=\al_1+\al_2,\;\;\bt_3=\al_2.$$

By Theorem \ref{nichols cy}, the algebra $R$ is CY if and only if
$$\chi_{_{\bt_2}}(g_1)\chi_{_{\bt_3}}(g_1)=(\chi_1\chi^2_2)(g_1)=1$$
and $$\chi_2(g_{_{\bt_1}})\chi_2(g_{_{\bt_2}})=\chi_2(g^2_1g_2)=1.$$
That is, $q_{_{11}}q_{_{12}}^2=q_{_{22}}q_{_{12}}^2=1$. 
By equation (\ref{q1}), we have
$q_{_{11}}^{-1}=q_{_{22}}^{-1}=q_{_{12}}q_{_{21}}. $

Now we conclude that the algebra $R$ is CY if and only if there is
some $q\in \kk^\times$, which is  not a root of unity, and satisfies the following relations:
$$q_{_{11}}=q_{_{22}}=q^2\;\t{ and }\;q_{_{12}}=q_{_{21}}=q^{-1}. $$
In other words, the braiding is of DJ-type. Then the algebra $R$ is an AS-regular algebra of type $A$ (see \cite{asch} for terminology). This coincides with Proposition 5.4 in \cite{br}.
\end{eg}

\begin{eg}
Let $R$ be a Nichols algebra of type $B_2$. That is, $R$ is
generated by $x_1$ and $x_2$  subjects the relations
$$\begin{array}{l}x_1^3x_2-q_{_{12}}x^2_1x_2x_1-q_{_{11}}q_{_{12}}x^2_1x_2x_1+q_{_{11}}q_{_{12}}^2x_1x_2x_1^2
\\-q^2_{_{11}}q_{_{12}}(x_1^2x_2x_1-q_{_{12}}x_1x_2x_1^2-q_{_{11}}q_{_{12}}x_1x_2x_1^2+q_{_{11}}q_{_{12}}^2x_2x_1^3)=0,\end{array}$$
and
$$x_2^2x_1-q_{_{21}}x_2x_1x_2-q_{_{22}}q_{_{21}}x_2x_1x_2+q_{_{22}}q_{_{21}}^2x_1x_2^2=0,$$
where $q_{ij}\in\kk$ for $1\se i,j\se 2$ and
$q_{_{12}}q_{_{21}}=q_{_{11}}^{-2}=q_{_{22}}^{-1}$. Applying a
similar argument, we obtain that $R$ is CY if and only if there is
some $q\in\kk^\times$, which is not a root of unity and satisfies the following:
$$q_{_{11}}=q, \;\;q_{_{12}}=q^{-1},\;\;q_{_{21}}=q^{-1} \t{ and }q_{_{22}}=q^{2}.$$
\end{eg}

\section{Relation between the Calabi-Yau property of pointed Hopf algebras and Nichols algebras}

We keep the notations as in the previous section. Let $(\mc{D},\lmd)$ be
a generic datum of finite Cartan type. In this section, we discuss the relation
between the  CY property of the algebra $U(\mc{D},\lmd)$ and that of
the corresponding Nichols algebra $\mc{B}(V)$. It turns out that if
one of them is CY, then the other one is not.

\begin{lem}\label{coeff}
For each $1\se k\se \tt$, we have
$$\prod_{i=1,i\neq j_k}^{p}\chi_{_{\bt_i}}(g_k)=(\prod_{i=1}^{j_k-1}\chi^{-1}_k(g_{_{\bt_i}}))( \prod_{i=j_k+1}^{p}\chi_{_{\bt_i}}(g_k)).$$
\end{lem}
\proof Let $\ome_0=s_{i_1}\cdots s_{i_p}$ be the fixed  reduced
decomposition of the longest element $\ome_0$ in the Weyl group. It
is clear that $\ome_0^{-1}$  is also of maximal length. By Lemma
3.11 in \cite{ka}, for each $1\se k\se \tt$, there exists $1\se t\se
p$, such that
$$s_{k}s_{i_1}\cdots s_{i_{t-1}}=s_{i_1}\cdots s_{i_t}.$$ That is,
$\ome_0=s_ks_{i_1}\cdots s_{i_{t-1}}s_{i_{t+1}}\cdots s_{i_p}$. Set
$$\bt_1'=\al_k, \;\;\bt_2'=s_k(\al_{i_1}),\cdots, \bt_p'=s_ks_{i_1}\cdots s_{i_{t-1}}s_{i_{t+1}}\cdots s_{i_{p-1}}(\al_{i_p}).$$
Applying a similar argument to the one  in the proof of Theorem \ref{nichols
cy}, we conclude that the rigid dualizing complex of the algebra
$R=\mc{B}(V)$ is isomorphic to $ _{\vph'}R[p]$. The algebra
automorphism $\vph'$ is defined by
$$\vph'(x_l)=(\prod_{i=1}^{j'_l-1}\chi^{-1}_l(g_{\bt_i'}))(
\prod_{i=j'_l+1}^{p}\chi_{\bt_i'}(g_l))x_l,$$ for each $1\se l\se
\tt$, where $j'_l$, $1\se l\se \tt$, are the integers such that
$\bt'_{j'_l}=\al_l$.  In particular, we have
$$\vph'(x_k)=( \prod_{i=2}^{p}\chi _{\bt_i'}(g_k))x_k.$$ The rigid dualizing complex is unique up to isomorphism, so $_{\vph'}R \cong {}_{\vph}R$ as $R$-$R$-bimodules, where $\vph$ is the algebra automorphism defined in Theorem \ref{nichols cy}. Since the graded algebra $R$ is connected, we have $\vph' =\vph $. In particular, $\vph'(x_k) =\vph(x_k)$, that is,  $$ \prod_{i=2}^{p}\chi _{\bt_i'}(g_k)=(\prod_{i=1}^{j_k-1}\chi^{-1}_k(g_{_{\bt_i}}))( \prod_{i=j_k+1}^{p}\chi _{_{\bt_i}}(g_k)).$$ Both $\bt_1,\cdots,\bt_p$ and $\bt_1',\cdots,\bt_p'$ are enumeration of positive roots. We have $\al_k=\bt_1'=\bt_{j_k}$. Therefore,
$$ \prod_{i=2}^{p}\chi _{\bt_i'}(g_k)=\prod_{i=1,i\neq j_k}^{p}\chi _{_{\bt_i}}(g_k).$$
It follows that $$(\prod_{i=1}^{j_k-1}\chi^{-1}_k(g_{_{\bt_i}}))(
\prod_{i=j_k+1}^{p}\chi _{_{\bt_i}}(g_k))=\prod_{i=1,i\neq
j_k}^{p}\chi _{_{\bt_i}}(g_k).$$\qed

\begin{prop}\label{prop cy-ni}
If $A=U(\mc{D},\lmd)$ is a CY algebra, then the rigid dualizing
complex of the Nichols algebra $R=\mc{B}(V)$ is isomorphic to $_\vph
R[p]$, where $\vph$ is defined by $\vph(x_k)=\chi^{-1}_k(g_k)x_k$,
for all $1\se k\se \tt$.
\end{prop}

\proof By Theorem \ref{nichols cy}, the rigid dualizing complex of
$R$ is isomorphic to $_\vph R[p]$, where $\vph$ is defined by
$$\vph(x_k)=(\prod_{i=1}^{j_k-1}\chi^{-1}_k(g_{_{\bt_i}}))(
\prod_{i=j_k+1}^{p}\chi _{_{\bt_i}}(g_k))x_k,$$ for all $1\se k\se \tt$.
If $A$ is a CY algebra, then
$\prod_{i=1}^{p}\chi_{_{\bt_i}}=\varepsilon$ by Theorem \ref{cy}.
Therefore, for $1\se k\se \tt$,
$$\begin{array}{ccl}(\prod_{i=1}^{j_k-1}\chi^{-1}_k(g_{_{\bt_i}}))( \prod_{i=j_k+1}^{p}\chi_{_{\bt_i}}(g_k))&=&\prod_{i=1,i\neq j_k}^{p}\chi _{_{\bt_i}}(g_k)\\
&=&\chi^{-1}_k(g_k),
\end{array}$$
where the first equation follows from Lemma \ref{coeff}. Now
$\vph(x_k)=\chi^{-1}_k(g_k)x_k$ for all $1\se k\se \tt$.  \qed

Note that $\chi_k(g_k)\neq 1$ for all $1\se k\se \tt$. So the algebra
$R=\mc{B}(V)$ is not CY, if $A=U(\mc{D},\lmd)$ is a CY algebra.

\begin{prop}\label{prop ni-cy}
If the Nichols algebra $R=\mc{B}(V)$ is a CY algebra, then the rigid
dualizing complex of $A=U(\mc{D},\lmd)$ is isomorphic to
$_{\psi}A[p+s]$, where $\psi$ is defined by $\psi(x_k)=x_k$ for all
$1\se k\se \tt$ and $\psi(g)=\prod_{i=1}^{p}\chi_{_{\bt_i}}(g)$ for
all $g\in \bgm$.
\end{prop}
\proof If the Nichols algebra $R$ is CY, then by Theorem \ref{nichols cy} and
Lemma \ref{coeff}, for each $1\se k\se \tt$, we have
$$\prod_{i=1,i\neq j_k}^{p}\chi_{_{\bt_i}}(g_k)=(\prod_{i=1}^{j_k-1}\chi^{-1}_k(g_{_{\bt_i}}))( \prod_{i=j_k+1}^{p}\chi_{_{\bt_i}}(g_k))=1.$$
Now the statement follows from
 Theorem \ref{cy}.\qed

With the assumption of Proposition \ref{prop ni-cy}, for all $1\se
k\se \tt$, we have
$$\psi(g_k)=\prod_{i=1}^{p}\chi_{_{\bt_i}}(g_k)=\chi_k(g_k)g_k\neq g_k.$$
Since the invertible elements of $A$ are in $\kk\bgm$ and $\bgm$ is
an abelian group,  $\psi$ can not be an inner automorphism. So the
algebra $A$ is not CY.

\begin{eg}
Let $R$ be the algebra in Example \ref{eg 1}. Assume that $\bgm
=\lan y_1,y_2\ran\cong \ZZ^2$, and  $g_i=y_i$, $i=1,2$. The
characters $\chi_1$ and $\chi_2$ are given by the following table:
$$\begin{tabular}{|c|c|c|}\hline
&$y_1$&$y_2$\\\hline $\chi_1$&$q^2$&$q^{-1}$\\\hline
$\chi_2$&$q^{-1}$&$q^2$\\\hline
\end{tabular}$$where $q$ is not a root of unity.

The algebra $R$ is a CY algebra. But the algebra $A=R\#\kk\bgm$ is
not. The rigid dualizing complex of $A$ is isomorphic to $_\psi
A[5]$, where $\psi$ is defined by $\psi(x_i)=x_i$ and
$\psi(y_i)=q^2y_i$ for $i=1,2$.
\end{eg}

\begin{eg}
Let $A$ be an algebra with generators $y_1^{\pm1}$, $y_2^{\pm1}$,
$x_1$ and $x_2$ subject to the relations
$$y_h^{\pm1}y_h^{\mp1}=1,\;\;1\se h,m\se 2,$$
$$y_1x_1=qx_1y_1,\;\;y_1x_2=q^{-1}x_2y_1,$$
$$y_2x_1=q^{\frac{k }{l }}x_1y_2,\;\;y_2x_2=q^{-\frac{k }{l }}x_2y_2,$$
$$x_1x_2-q^{-k}x_2x_1=1-y_1^{k }y_2^{l },$$
where $k,l \in \ZZ^+$ and $q\in \kk$ is not a root of unity.

By Theorem \ref{cy}, the algebra $A$ is a CY algebra of dimension 4. Let $R$ be
the corresponding Nichols algebra of $A$. The rigid dualizing complex
of $R$ is isomorphic to $_\vph R[2]$, where $\vph$ is defined by
$\vph (x_1)=q^{-k}x_1$ and $\vph (x_2)=q^k x_2$.
\end{eg}

\section{Classification of Calabi-Yau pointed Hopf algebra $U(\mc{D},\lmd)$ of lower dimensions}
In this section, we assume that $\kk=\mathbb{C}$. We shall classify CY pointed Hopf algebras
$U(\mc{D},\lmd)$ of dimension less than 5, where $(\mc{D},\lmd)$ is a
generic datum of finite Cartan type. In a generic datum $(\bgm,(a_{ij}),(q_{_I}),
(g_i),(\chi_i),(\lmd_{ij}))$ of finite Cartan type, $(q_{_I})$ are determined by $(\chi_i)$ and $(g_i)$. In the following,  we will omit $(q_{_I})$ for simplicity.

Let $(\mc{D},\lmd)=(\bgm,(a_{ij}),(g_i),(\chi_i),(\lmd_{ij}))$  be  a generic datum of
finite Cartan type. Then $\chi_i(g_i)$ are not roots of unity
for $1\se i\se \tt$. Hence, in the classification, we exclude the
case where the group is trivial. If the group $\bgm$ in a datum
$(\mc{D},\lmd)=(\bgm,(a_{ij}),
(g_i),(\chi_i),(\lmd_{ij}))$ is trivial, then the algebra
$U(\mc{D},0)$  (in this case, $U(\mc{D},0)$ has no non-trivial
lifting) is the universal enveloping algebra $U(\mathfrak{g})$,
where the Lie algebra $\mathfrak{g}$ is generated by $x_i$, $1\se
i\se \tt$, subject to the relations
$$(\ad x_i)^{1-a_{ij}}x_j=0,\;\;1\se i,j\se \tt,\;\;i\neq j.$$
We have $\t{tr}(\ad x)=0$ for all $x\in \mathfrak{g}$. Therefore,
$U(\mathfrak{g})$ is CY by \cite[Lemma 4.1]{hvz}.  We list those of
dimension less than 5 in the following table.

\begin{center}
\begin{longtable}{|c|c|c|c|c|}
\hline &CY&  & \multicolumn{2}{c}{Lie algebra} \vline   \\
\cline{4-5} Case&dimension &Cartan matrix&bases&relations\\ \hline
1&1&$A_1$&$x$&\\\hline 2&2&$A_1\times A_1$&$x,y$&abelian Lie
algebra\\\hline

3&3&$A_1\times A_1\times A_1$&$x,y,z$& abelian Lie algebra\\ \hline

4&3&$A_2$&$x,y,z$&$[x,y]=z,[x,z]=[y,z]=0$\\ \hline

5&4&$A_1\times\cdots\times A_1$&$x,y,z,w$& abelian Lie algebra\\
\hline
6&4&$A_1\times A_2$&$x,y,z,w$& $[x,y]=z,[x,z]=[y,z]=0$\\
&&&&$[x,w]=[y,w]=[z,w]=0$\\\hline
7&4&$B_2$ &$x,y,z,w$& $[x,y]=z,[x,z]=w$,\\
&&&& $[x,w]=[y,z]$ \\
&&&& $=[y,w]=[z,w]=0$ \\\hline
\end{longtable}
\end{center}

\begin{rem}
The Lie algebra in Case 4 is the Heisenberg algebra. In \cite{hvz},
the authors classified those 3-dimensional Lie algebras whose
universal enveloping algebras are CY algebras. Beside the algebras
in Case 3 and Case 4, the other two Lie algebras are
\begin{itemize}
\item  The 3-dimensional simple Lie algebra $\mk{sl}_2$;
\item  The Lie algebra $\mk{g}$, where $\mk{g}$ has a basis
\{x,y,z\} such that $[x,y]=y$, $[x,z]=-z$ and $[y,z]=0$.
\end{itemize}
\end{rem}

Now let $$(\mc{D},\lmd)=(\bgm,(a_{ij})_{1\se i,j\se \tt},
(g_i)_{1\se i\se \tt},(\chi_i)_{1\se i\se \tt},(\lmd_{ij})_{1\se
i<j\se \tt,i\nsim j} )$$  and
$$(\mc{D}',\lmd')=(\bgm',
(a'_{ij})_{1\se i,j\se \tt'},   (g'_i)_{1\se
i\se \tt',},(\chi'_i)_{1\se i\se \tt'},(\lmd'_{ij})_{1\se i<j\se
\tt',i\nsim j})$$ be two generic data of finite Cartan type for groups $\bgm$
and $\bgm'$, where $\bgm$ and $\bgm'$ are both free abelian groups
of finite rank.

The data $(\mc{D},\lmd)$ and $(\mc{D}',\lmd')$ are said to be
\textit{isomorphic} if $\tt=\tt'$, and if there exist a group
isomorphism $\vph:\bgm\ra \bgm'$, a permutation $\sigma\in
\mathbb{S}_\tt$, and elements $0\neq \al_i\in \kk$, for all $1\se
i\se \tt$  subject to the following relations:
\begin{enumerate}
\item[(i)] $\vph(g_i)=g'_{\sigma(i)}$,  for all $1\se i\se \tt$.
\item[(ii)] $\chi_i=\chi'_{\sigma(i)}\vph$, for all $1\se i\se \tt$.
\item[(iii)] $\lmd_{ij}=\begin{cases}\al_i\al_j\lmd'_{\sigma(i)\sigma(j)}, & \t{if } \sigma(i)<\sigma(j)\\-\al_i\al_j\chi_j(g_i)\lmd'_{\sigma(j)\sigma(i)}, & \t{if }\sigma(i)>\sigma(j)\end{cases}$,
   \vspace{2mm}

    for all  $1\se i<j\se \tt$ and $i\nsim j$.
\end{enumerate}
In this case the triple $(\vph,\sigma, (\al_i))$ is called an
\textit{isomorphism} from $(\mc{D},\lmd)$ to $(\mc{D}',\lmd')$.

If $(\mc{D},\lmd)$ and $(\mc{D}',\lmd)$ are isomorphic, then we can
deduce that $a_{ij}=a'_{\sg(i)\sg(j)}$ for all $1\se i,j\se \tt$
\cite{as3}.

The following corollary can be immediately obtained from the
definition of isomorphic data.
\begin{cor}\label{isorem}
Let ($\mc{D}$,0) be a generic datum of finite Cartan type formed by
$(\bgm,(a_{ij}),
(g_i),(\chi_i),0)$.  Assume that
$\vph:\bgm\ra \bgm'$ is  a group isomorphism and $\sigma$ is a
permutation in $\mathbb{S}_\tt$. Then $(\mc{D}, 0)$ is isomorphic to
 $(\mc{D}', 0)$, where $\mc{D}'$ is formed by
 $(\bgm',(a_{\sg^{-1}(i)\sg^{-1}(j)}), (\vph(g_{\sg^{-1}(i)})) ,(\chi_{\sg^{-1}(i)}\vph^{-1}))$.
\end{cor}

Let $(\mc{D},\lmd)$ be a generic datum of finite Cartan type. Following from \cite{as3}, the pointed Hopf algebra
$U(\mc{D},\lmd)$ is uniquely determined by datum $(\mc{D},\lmd)$.
Let $Isom((\mc{D},\lmd),(\mc{D}',\lmd'))$ be the set of all
isomorphisms from $(\mc{D},\lmd)$ to $(\mc{D}',\lmd')$. Let $A,B$ be
two Hopf algebras, we denote by $Isom(A,B)$  the set of all Hopf
algebra isomorphisms from $A$ to $B$.

\begin{lem}\label{iso}\cite[Thm. 4.5]{as3}
Let $(\mc{D},\lmd)$ and $(\mc{D}',\lmd')$ be two generic data of
finite Cartan type. Then the Hopf algebras $U(\mc{D},\lmd)$ and
$U(\mc{D}',\lmd')$  are isomorphic if and only if $(\mc{D},\lmd)$ is
isomorphic to $(\mc{D}',\lmd')$. More precisely, let
$x_1,\cdots,x_\tt$ (resp. $x'_1,\cdots,x'_\tt$) be the simple root
vectors in $U(\mc{D},\lmd)$ (resp. $U(\mc{D'},\lmd')$), and let
$g_1,\cdots,g_\tt$ (resp. $g'_1,\cdots,g'_\tt$) be the group-like
elements in $\mc{D}$ (resp. $\mc{D}'$).  Then the map
$$Isom(U(\mc{D},\lmd), U(\mc{D'},\lmd')) \ra  Isom((\mc{D},\lmd),(\mc{D'},\lmd')),$$
given by $\phi\mapsto(\vph,\sigma,(\al_i))$, where
$\vph(g)=\phi(g)$, $\vph(g_i)=g'_{\sigma(i)}$,
$\phi(x_i)=\al_ix'_{\sigma(i)}$, for all $g\in \bgm$, $1\se i\se
\tt$, is bijective.
\end{lem}

The following lemma is well-known.
\begin{lem}\label{poly}
If $\bgm$ is a free abelian group of rank $s$, then the algebra $\kk
\bgm$ is a CY algebra of dimension $s$.\end{lem}


If $\bgm$ is a free abelian group of finite rank, we denote by
$|\bgm|$ the rank of $\bgm$.

\begin{prop}\label{12}
Let $A$ be the algebra $U(\mc{D},\lmd)$, where $(\mc{D},\lmd)$ is a
generic datum of finite Cartan type for a group $\bgm$. Then
\begin{enumerate} \item[(a)] $A$ is
CY of dimension 1 if and only if $A=\kk \ZZ$.
\item[(b)] $A$ is CY of dimension 2 if and only if $A=\kk \bgm$, where $\bgm$ is a free abelian group
of rank 2.

\end{enumerate}
\end{prop}
\proof (a) is clear.

(b) It is sufficient to show that if $A$ is CY of dimension 2, then
$A$ is the group algebra of a free abelian group of rank 2.  By
Theorem \ref{as}, if the global dimension of $A$ is 2. Then the
following possibilities arise:
\begin{enumerate}
\item[(i)] $|\bgm|=2$, $A=\kk \bgm$ is the group algebra of a free abelian group of rank 2;
\item[(ii)] $|\bgm|=1$ and the Cartan matrix of $A$ is of type $A_1$.
\end{enumerate}

Let $A$ be a pointed Hopf  algebra  of type (ii) and  let the datum
$$(\mc{D},\lmd)=(\bgm,  (g_i),(\chi_i),(a_{ij}) , (\lmd_{ij}) )$$ be
as follows
\begin{itemize}
\item $\bgm=\langle y_1\rangle\cong \ZZ$;
\item $g_1=y_1^k$, for some $k\in \ZZ$;
\item $\chi_1\in \widehat{\bgm}$ is defined by $\chi_1(y_1)=q$, where $q$ is not a root of unity;
\item The Cartan matrix is of type $A_1$;
\item $\lmd=0$.
\end{itemize}
Observe that in this case, the linking parameter must be 0. In
addition, there is only one root vector, that is, the simple root
vector $x_1$. Since $q\neq 1$, we have $\chi_1\neq \varepsilon$. So
the algebra $A$ is not CY by Theorem \ref{cy}.

Therefore, if $A$ is CY, then $A$ is of type (i). Hence, the
classification is complete. \qed

\begin{prop}\label{3}
Let $A$ be the algebra $U(\mc{D},\lmd)$, where $(\mc{D},\lmd)$ is a
generic datum of finite Cartan type for a group $\bgm$. If $A$ is CY
of dimension 3, then the group $\bgm$ and the Cartan matrix
$(a_{ij})$ are given by one of the following 2
 cases.
\begin{center}
\begin{longtable}{ccc}
\hline Case&{$|\bgm|$}&{Cartan matrix}  \\ \hline
1&3 &Trivial\\
2&1&$A_1\times A_1$\\
\hline
\end{longtable}
\end{center}

The non-isomorphic classes of CY algebras in each case are given as
follows.

Case 1: The group algebra of a free abelian group of rank 3.

Case 2: \begin{enumerate}\item[(I)]The datum $(\mc{D},\lmd)=(\bgm,
(g_1,g_2),(\chi_1,\chi_2),(a_{ij})_{1\se i,j\se 2},\lmd_{12})$ is
given as follows:
\begin{itemize}
\item $\bgm=\lan y_1\ran\cong \ZZ$;
\item $g_1=g_2=y_1^{k}$ for some $k\in \ZZ^+$;\item  $\chi_1(y_1)=q$, where $q\in\kk$  is not a root of
unity and $0<|q|<1$, and $\chi_2=\chi_1^{-1}$;
\item $(a_{ij})_{1\se
i,j\se 2}$ is the Cartan matrix  of type $A_1\times A_1$;
\item $\lmd_{12}=0$.
\end{itemize}

\item[(II)]The datum
$(\mc{D},\lmd)=(\bgm, (g_1,g_2),(\chi_1,\chi_2),(a_{ij})_{1\se
i,j\se 2},\lmd_{12})$ is given as follows:
\begin{itemize}
\item $\bgm=\lan y_1\ran\cong \ZZ$;
\item $g_1=g_2=y_1^{k}$ for some $k\in \ZZ^+$;
\item  $\chi_1(y_1)=q$, where $q\in\kk$  is not a root of
unity and $0<|q|<1$, and $\chi_2=\chi_1^{-1}$;
\item $(a_{ij})_{1\se i,j\se 2}$ is the Cartan matrix of type $A_1\times A_1$;
\item $\lmd_{12}=1$.
\end{itemize}
\end{enumerate}
\end{prop}
\proof  By Remark \ref{lifting}, it is sufficient to  discuss the
graded case and consider the non-trivial liftings. We first show
that the algebras listed in the proposition are all CY. Case 1
follows from Lemma \ref{poly}. Now we discuss Case 2. The root
system of the Cartan matrix of type $A_1\times A_1$ has two simple
roots, say $\al_1$ and $\al_2$. They are also the positive roots.
First we have $\chi_1\chi_2=\varepsilon$. Since
$\mc{S}^2_A(x_i)=\chi_i(g^{-1}_i)x_i$, $i=1,2$, $g_1=g_2=y_1^{k}$,
we have $\mc{S}^2_A(x_i)={y_1}^{-k}x_i{y_1}^k$ for $i=1,2$. It is
easy to see that $\mc{S}_A^2(y_1)=y_1$. It follows that $\mc{S}^2_A$
is an inner automorphism. Thus the algebras in Case 2 are CY by
Theorem \ref{cy}.

Now we show that the classification is complete.

If $A$ is  of   global dimension 3, then the following possibilities
for the group $\bgm$ and the Cartan matrix $(a_{ij})$ arise:
\begin{enumerate}
\item[(i)] $|\bgm|=3$, $A$ is the group algebra of a free abelian group of rank 3.
\item[(ii)] $|\bgm|=2$ and the Cartan matrix of $A$ is of type $A_1$.
\item[(iii)] $|\bgm|=1$ and the Cartan matrix of $A$ is of type $A_1\times A_1$.
\end{enumerate}
Similar to the case of global dimension 2, $A$ can not be CY if $A$
is of type (ii).

Now, let  $A$ be a CY graded algebra of type (iii). In this case, we
have $\chi_2(g_1)\chi_1(g_2)=1$ (cf. equation (\ref{q})). In
addition, we have $\chi_1\chi_2=\varepsilon$ by Theorem \ref{cy}. It
follows that $1=\chi_2(g_1)\chi_1(g_2)
=\chi^{-1}_1(g_1)\chi_1(g_2)$. Let $\bgm=\lan y_1\ran$ and
assume that $g_1=y_1^k$, $g_2=y_1^l$ for some $k,l\in \mathbb{Z}$. Then
$\chi_1(y_1^{l-k})=1$. Since $\chi_1(y_1)$ is not a root of unity,
we have $k=l$, that is, $g_1=g_2=y_1^k$. Therefore, $A\cong
U(\mc{D},0)$, where the datum $\mc{D}$ is given by
\begin{itemize}
\item $\bgm=\lan y_1\ran\cong \ZZ$;
\item $g_1=g_2=y_1^{k}$, for some $k\in \ZZ$;\item  $\chi_1(y_1)=q$, where $q\in\kk$  is not a root of unity, and $\chi_2=\chi_1^{-1}$;
\item $(a_{ij})_{1\se i,j\se 2}$ is the Cartan matrix of type $A_1\times
A_1$.
\end{itemize}

Let  $ \mc{D}' $ be another datum given by
\begin{itemize}
\item $\bgm'=\lan y_1'\ran\cong \ZZ$;
\item $g_1'=g_2'=y_1'^{k'}$, for some $k'\in \ZZ$;\item  $\chi_1'(y_1')=q'$, where $q'\in\kk$  is not a root of unity, and $\chi_2'=\chi_1'^{-1}$;
\item $(a'_{ij})_{1\se i,j\se 2}$ is the Cartan matrix of type $A_1\times
A_1$.
\end{itemize}
Assume that $(\mc{D}', 0)$  is isomorphic to  $(\mc{D},0)$ via an
isomorphism $(\vph,\sigma,(\al_i))$. Then  $\vph$ is a group
automorphism such that $\vph(y_1)=y_1'$ or $\vph(y_1)=y_1'^{-1}$.
Since $\sg\in \mathbb{S}_2$, we have $\sg=\id$ or $\sg=(12)$. From
an easy computation, there are four possibilities for $k'$ and $q'$,
\begin{itemize}
\item $k'=k$ and $q'=q$;
\item $k'=-k$ and $q'=q$;
\item $k'=k$ and $q'=q^{-1}$;
\item $k'=-k$ and $q'=q^{-1}$.
\end{itemize}
This shows that $A= U(\mc{D},0)$ is isomorphic to  an algebra in (I)
of Case 2. In addition,  every pair $(k,q)\in \ZZ^+\times \kk$, such
that $0<|q|< 1$ determines a non-isomorphic algebra in (I) of Case 2.
Each algebra in (I) of Case 2 has only one non-trivial lifting,
which is isomorphic to an algebra in (II).

Thus we have completed the classification. \qed

We list all CY  Hopf algebras  $U(\mc{D},\lmd)$ of dimension 3 in
terms of  generators and relations in the following table. Note that
in each case  $q$ is not a root of unity.

\begin{center}
CY algebras of dimension 3
\begin{longtable}{ccc}
\hline Case &Generators&Relations \\ \hline Case 1& $y_h,y_h^{-1}$ &
$y_h^{\pm1}y_m^{\pm1}=y_m^{\pm1}y_h^{\pm1}$\\
&$1\se h\se 3$&$y_h^{\pm1}y_h^{\mp1}=1$\\
&&$1\se h,m\se 3$\\\hline

Case 2 (I)&$y_1^{\pm1},x_1,x_2$&$y_1y_1^{-1}=y_1^{-1}y_1=1$\\
&&$y_1x_1=qx_1y_1$\\
&&$y_1x_2=q^{-1}x_2y_1$, $0<|q|<1$\\
&&$x_1x_2-q^{-k}x_2x_1=0$, $k\in \ZZ^+$\\\hline

Case 2 (II)&$y_1^{\pm1},x_1,x_2$&$y_1y_1^{-1}=y_1^{-1}y_1=1$\\
&&$y_1x_1=qx_1y_1$\\
&&$y_1x_2=q^{-1}x_2y_1$, $0<|q|<1$\\
&&$x_1x_2-q^{-k}x_2x_1=(1-y_1^{2k})$, $k\in \ZZ^+$\\
\hline


\end{longtable}
\end{center}

\begin{prop}\label{4}
Let $A$ be the algebra $U(\mc{D},\lmd)$, where $(\mc{D},\lmd)$ is a
generic datum of finite Cartan type for a group $\bgm$. If  $A$ is
CY of dimension 4, then the group $\bgm$ and the Cartan matrix
$(a_{ij})$ are given by one of the following 2 cases.
\begin{center}
\begin{longtable}{ccc}
\hline Case&{$|\bgm|$}&{Cartan matrix} \\ \hline
1&4&Trivial\\
2&2&$A_1\times A_1$\\
\hline
\end{longtable}
\end{center}

In  each case, the non-isomorphic classes of CY algebras are given
as follows.

Case 1: The group algebra of a free abelian group of rank 4.

Case 2: \begin{enumerate}\item[(I)] The datum $(\mc{D},\lmd)=(\bgm,
(g_1,g_2),(\chi_1,\chi_2),(a_{ij})_{1\se i,j\se 2},\lmd_{12})$ is
given by
\begin{itemize}
\item $\bgm=\lan y_1,y_2 \ran\cong \ZZ^2$;
\item $g_1=g_2=y_1^{k}$ for some $k\in \ZZ^+$;
\item  \begin{itemize} \item $\chi_1(y_1)=q_{_1}$, $\chi_1(y_2)=q_{_2}$, where  $q_{_1},q_{_2}\in\kk$ satisfy that
$0<|q_{_1}| <1$ and $q_{_1}$ is not a root of
unity,\item$\chi_2=\chi_1^{-1}$;
\end{itemize}
\item $(a_{ij})_{1\se i,j\se 2}$ is the Cartan matrix of type $A_1\times A_1$;
\item $\lmd_{12}=0$.
\end{itemize}
\item[(II)] The datum
$(\mc{D},\lmd)=(\bgm, (g_1,g_2),(\chi_1,\chi_2),(a_{ij})_{1\se
i,j\se 2},\lmd_{12})$ is given by
\begin{itemize}
\item $\bgm=\lan y_1,y_2 \ran\cong \ZZ^2$;
\item $g_1=g_2=y_1^{k}$ for some $k\in \ZZ^+$;
\item   \begin{itemize} \item $\chi_1(y_1)=q_{_1}$, $\chi_1(y_2)=q_{_2}$, where  $q_{_1},q_{_2}\in\kk$ satisfy that
$0<|q_{_1}|<1$ and $q_{_1}$ is not a root of
unity,\item$\chi_2=\chi_1^{-1}$;
\end{itemize}
\item $(a_{ij})_{1\se i,j\se 2}$ is the Cartan matrix of type $A_1\times A_1$;
\item $\lmd_{12}=1$.
\end{itemize}
\end{enumerate}

Let $A$ and $B$ be two algebras in Case (I) (or (II)) defined by
triples $(k, q_{_1}, q_{_2})$ and $(k', q_{_1}', q_{_2}')$ respectively. They
are isomorphic if and only if $k=k'$, $q_{_1}=q_{_1}'$ and there is some
integer $b$, such that $q_{_2}'=q_{_1}^bq_{_2}$ or $q_{_2}'=q_{_1}^bq_{_2}^{-1}$.
\begin{enumerate}

\item[(III)] The datum
$(\mc{D},\lmd)=(\bgm, (g_1,g_2),(\chi_1,\chi_2),(a_{ij})_{1\se
i,j\se 2},\lmd_{12})$ is given by
\begin{itemize}
\item $\bgm=\lan y_1,y_2\ran\cong\ZZ^2$;
\item $g_1=y_1^{k}$, $g_2= y_2^{l }$ for some $k, l \in \ZZ^+$; \item \begin{itemize}\item$\chi_1(y_1)=q$, $\chi_1(y_2)=q^{\frac{k }{l }}$, where $q\in\kk$  is not a root of unity and $0<|q|<1$, \item
$\chi_2=\chi_1^{-1}$;\end{itemize}
\item $(a_{ij})_{1\se i,j\se 2}$ is the Cartan matrix of type $A_1\times A_1$;
\item $\lmd_{12}=0$.
\end{itemize}
\item[(IV)] The datum
$(\mc{D},\lmd)=(\bgm, (g_1,g_2),(\chi_1,\chi_2),(a_{ij})_{1\se
i,j\se 2},\lmd_{12})$ is given by
\begin{itemize}
\item $\bgm=\lan y_1,y_2\ran\cong\ZZ^2$;
\item $g_1=y_1^{k}$, $g_2= y_2^{l }$ for some $k, l \in \ZZ^+$; \item \begin{itemize}\item$\chi_1(y_1)=q$, $\chi_1(y_2)=q^{\frac{k }{l }}$, where $q\in\kk$  is not a root of unity and  $0<|q|<1$,  \item
$\chi_2=\chi_1^{-1}$;\end{itemize}\item $(a_{ij})_{1\se i,j\se 2}$
is the Cartan matrix of type $A_1\times A_1$;
\item $\lmd_{12}=1$.
\end{itemize}

\item[(V)] The datum
$(\mc{D},\lmd)=(\bgm, (g_1,g_2),(\chi_1,\chi_2),(a_{ij})_{1\se
i,j\se 2},\lmd_{12})$ is given by
\begin{itemize}
\item $\bgm=\lan y_1,y_2\ran\cong\ZZ^2$;
\item $g_1=y_1^{k}$, $g_2=y_1^{l_1}y_2^{l_2}$ for some $k, ,l_1,l_2\in \ZZ^+$, $k\neq l_1$, $0< l_1<l_2$; \item\begin{itemize}\item $\chi_1(y_1)=q$, $\chi_1(y_2)=q^{\frac{k-l_1}{l_2}}$, where $q\in\kk$  is not a root of unity and $0<|q|<1$,
\item $\chi_2=\chi_1^{-1}$;\end{itemize}
\item $(a_{ij})_{1\se i,j\se 2}$ is the Cartan matrix of type $A_1\times A_1$;
\item $\lmd_{12}=0$.
\end{itemize}
\item[(VI)] The datum
$(\mc{D},\lmd)=(\bgm, (g_1,g_2),(\chi_1,\chi_2),(a_{ij})_{1\se
i,j\se 2},\lmd_{12})$ is given by
\begin{itemize}
\item $\bgm=\lan y_1,y_2\ran\cong\ZZ^2$;
\item $g_1=y_1^{k}$, $g_2=y_1^{l_1}y_2^{l_2}$ for some $k,l_1, l_2\in \ZZ^+$, $k\neq l_1$ and   $0<l_1<l_2$; \item\begin{itemize}\item $\chi_1(y_1)=q$, $\chi_1(y_2)=q^{\frac{k-l_1}{l_2}}$, where $q\in\kk$  is not a root of unity and  $0<|q|<1$,
\item $\chi_2=\chi_1^{-1}$;\end{itemize}\item $(a_{ij})_{1\se i,j\se 2}$ is the Cartan matrix of type $A_1\times A_1$;
\item $\lmd_{12}=1$.
\end{itemize}

\end{enumerate}
\end{prop}

\proof   We first show that the algebras listed in the proposition
are all CY. That the algebra in Case 1  is a CY algebra follows from
Lemma \ref{poly}. In  Case 2,  we have $\chi_1\chi_2=\varepsilon$
and $\mc{S}_A^2$ is an inner automorphism in each subcase. Indeed,
$\mc{S}^2_A(x_i)=g_1^{-1}x_ig_1$ and
$\mc{S}_A^2(y_i)=g_1^{-1}y_ig_1=y_i$, $i=1,2$. Thus the algebras in
Case 2 are CY by Theorem \ref{cy}.

Now we show that the classification is complete and the algebras on
the list are non-isomorphic to  each other.

If $A$ is  of global dimension 4, then the  group $\bgm$ and the
Cartan matrix $(a_{ij})$ must be one of the following types:
\begin{enumerate}
\item[(i)] $|\bgm|=4$ and $A$ is the group algebra of a free abelian group of rank 4.
\item[(ii)] $|\bgm|=3$ and the Cartan matrix of $A$ is of type $A_1$.
\item[(iii)] $|\bgm|=2$ and the Cartan matrix of $A$ is of type $A_1\times A_1$.
\item[(iv)] $|\bgm|=1$ and the Cartan matrix of $A$ is of type $A_1\times A_1\times A_1$.
\item[(v)] $|\bgm|=1$ and the Cartan matrix of $A$ is of type $A_2$.
\end{enumerate}

Let $A$ be a CY algebra of dimension 4.  Similar to the case of
global dimension 2, $A$ cannot be  of type (ii). We claim that $A$
cannot be  of type (iv) and (v) either.

Assume that $A$ is of type (iv), put $\bgm=\lan y_1\ran$, $g_i=y_1^{m_i}$
for some $0\neq m_i\in \ZZ$ and $\chi_i(y_1)=q_i$ for some
$q_i\in\kk$, $1\se i\se 3$. Then $q_{ij}=q_j^{m_i}$, for $1\se
i,j\se 3$. Because each $q_{ii}$ is not a root of unity, each $q_i$
is not a root of unity either.    Since $q_{ij}q_{ji}=1$, we have
$$q_{_1}^{m_2}q_{_2}^{m_1}=1,\;\;q_{_1}^{m_3}q_{_3}^{m_1}=1,\;\;q_{_2}^{m_3}q_{_3}^{m_2}=1.$$ Then $q_{_1}^{2{m_2m_3}}=1$. But
$q_{_1}$ is not a root of unity. So $A$ can not be  of   type (4).

In case (v), there are 3 positive roots in the root system. They are
$\al_1$, $\al_2$ and $\al_1+\al_2$, where $\al_1$ and $\al_2$ are
the simple roots. If $A$ is CY, then $\chi_1^2\chi_2^2=\varepsilon$
by Theorem \ref{cy}. So  we have $q_{_{11}}^2q_{_{21}}^2=1$ and
$q_{_{12}}^2q_{_{22}}^2=1$. However, $q_{_{21}}q_{_{12}}=q_{_{11}}^{-1}$ (equation
(\ref{q1})). Thus $q^2_{_{22}}=1$. But $q_{_{22}}$ is not a root of
unity. So $A$ cannot be of type (v) either.

Now to show that the classification is complete, we only need to
show that if $A$ is a  CY pointed Hopf algebra of type (iii), then $A$
is isomorphic to an algebra in   Case 2. Each algebra in (I), (III)
and (V) of Case 2 has only one non-trivial lifting, which is
isomorphic to an algebra in (II), (IV) and (VI) respectively. By
Remark \ref{lifting}, it suffices to  show that if $A$ is a graded
CY pointed Hopf algebra of type (iii), then $A$ is isomorphic to an
algebra in (I), (III) and (V) of Case 2.

Let $\bgm=\lan y_1,y_2\ran$ be  a free abelian group of rank 2. We
write $\chi_1(y_1)=q_{_1}$, $\chi_1(y_2)=q_{_2}$ and
$g_1=y_1^{k_1}y_2^{k_2}$, $g_2=y_1^{l_1}y_2^{l_2}$,  where
$\chi_1(g_1)=q_{_1}^{k_1}q_{_2}^{k_2}$ is  not a root of unity, and
 $k_1,k_2,l_1,l_2\in
\ZZ$. Following Theorem \ref{cy}, we have
$\chi_1\chi_2=\varepsilon$. So $q_{_{21}}=q_{_1}^{l_1}q_{_2}^{l_2}$
and $q_{_{12}}=q_{_1}^{-k_1}q_{_2}^{-k_2}$. We also have
$q_{_{12}}q_{_{21}}=1$ (equation (\ref{q1})). Thus
$q_{_1}^{l_1-k_1}q_{_2}^{l_2-k_2}=1$. Therefore,
 $A\cong U(\mc{D},0)$, where the datum
$\mc{D}$ is formed  by
\begin{itemize}\item $\bgm=\lan y_1,y_2\ran\cong\ZZ^2$; \item $(a_{ij})$ is the Cartan matrix of type
$A_1\times A_1$;\item $g_1=y_1^{k_1}y_2^{k_2}$,
$g_2=y_1^{l_1}y_2^{l_2}$, $k_1,k_2,l_1,l_2\in \ZZ$;\item
$\chi_1(y_1)=q_{_1}$, $\chi_1(y_2)=q_{_2}$, where
$\chi_1(g_1)=q_{_1}^{k_1}q_{_2}^{k_2}$ is  not a root of unity and
$q_{_1}^{l_1-k_1}q_{_2}^{l_2-k_2}=1$, and
 $\chi_2=\chi_1^{-1}$.\end{itemize}

In the above datum $\mc{D}$, we may  assume that $k_1>0$ and
$k_2=0$. Then $q_{_1}$ is not a root of unity.  We show that there
is a group isomorphism $\vph:\bgm\ra \bgm'$, where $\bgm'=\lan
y'_1,y'_2\ran $ is also a free abelian group of rank 2, such that
$\vph(y_1^{k_1}y_2^{k_2})=y_1'^k$ and $k>0$.

The integers $k_1$ and $k_2$ can not be both equal to 0. If $k_2=0$
and $k_1>0$, then it is done. If $k_2=0$ and $k_1<0$, then
$\vph(y_1)=y_1'^{-1}$ and $\vph(y_2)=y_2'^{-1}$ defines a desired
isomorphism.

Similarly, we can obtain a desired isomorphism when $k_1=0$ and
$k_2\neq 0$.

If $k_1,k_2\neq 0$, then there are some $k, \bar{k}_1,\bar{k}_2\in
\ZZ$, such that $k_1=\bar{k}_1k$, $k_2=\bar{k}_2k$,   $k>0$ and
$(\bar{k}_1,\bar{k}_2)=1$, that is, $\bar{k}_1$ and $\bar{k}_2$ have
no common divisors. We can find  integers $a,b$ such that
$a\bar{k}_1+b\bar{k}_2=1$. Let $\vph:\bgm\ra \bgm'$ be the group
isomorphism defined by $\vph(y_1)={y_1'}^a{y_2'}^{-\bar{k}_2}$ and
$\vph(y_2)={y_1'}^b{y_2'}^{\bar{k}_1}$. Then
$\vph(y_1^{k_1}y_2^{k_2})={y_1'}^k$ and $k>0$. In conclusion, we
have proved the claim.

If $l_2=0$, then we have $q_{_1}^{l_1-k_1}=1$. Since $q_{_1}$ is not a
root of unity, we have $l_1=k_1$.  Applying a similar argument to
the one in Case 2 of Proposition \ref{3}, we find that $A$ is
isomorphic to  an algebra in (I) of Case 2.

Next, we consider the case when $l_2\neq 0$. In case  $l_1=0$, like
what we did for $k_1$ and $k_2$, we may assume that  $l_2>0$. If
$0<|q_{_1}|<1$, then $A$ is isomorphic to an algebra in (III) of Case 2.
Otherwise, the datum $(\mc{D},0)$  is isomorphic to  the datum given
by
\begin{itemize}
\item $\bgm=\lan y_1,y_2\ran\cong \ZZ^2$;
\item $g_1'=y_1^{l_2}$,
$g_2'=y_2^{k_1}$, $k_1,l_2\in \ZZ^+$;
\item $\chi_1'(y_1)=q_{_1}^{-\frac{k_1}{l_2}}$,
$\chi_1'(y_2)=q_{_1}^{-1}$,
 $\chi_2'=\chi_1'^{-1}$.
\item $(a_{ij})$ is the Cartan matrix of type
$A_1\times A_1$;
\item $\lmd_{12}=0$\end{itemize}
via the isomorphism $(\vph,(12),\al_1=\al_2=1)$, where $\vph$ is the
algebra automorphism defined by $\vph(y_1)=y_2$ and $\vph(y_2)=y_1$.
So $A$ is isomorphic to an algebra in (III) of Case 2 as well.

If $l_1\neq 0$ and $l_2>0$, then there is an integer $c$, such that
$0\se l_1+cl_2<l_2$.  Let $\bgm'=\lan y'_1,y'_2\ran $ be a free
abelian group of rank 2, and $\vph:\bgm\ra \bgm'$  the group
isomorphism defined by $\vph(y_1)=y_1'$ and
$\vph(y_2)=y_1'^{c}y_2'$. Then $\vph(y_1^{ k_1})={y_1'}^{k_1}$ and
$\vph(y_1^{l_1}y_2^{l_2})={y_1'}^{l_1+cl_2}{y_2'}^{l_2}$.

If $l_1\neq 0$ and $l_2<0$, then there are integers $\bar{l}_1$,
$\bar{l}_2$, such that $l_1=\bar{l}_1l$, $l_2=\bar{l}_2l$, $l>0$ and
$(\bar{l}_1,\bar{l}_2)=1$. So $\bar{l}_2<0$.  We can find integers
$a,b$ such that $a\bar{l}_1+b\bar{l}_2=1$. Since for any integer
$d$,
$(a+d\bar{l}_2)\bar{l}_1+(b-d\bar{l}_1)\bar{l}_2=a\bar{l}_1+b\bar{l}_2=1$,
we may assume that $0\se a<-\bar{l}_2$. Let $\bgm'=\lan
y'_1,y'_2\ran $ be  a free abelian group of rank 2, and
$\vph:\bgm\ra \bgm'$ be the group isomorphism defined by
$\vph(y_1)={y_1'}^a{y_2'}^{-\bar{l}_2}$ and
$\vph(y_2)={y_1'}^b{y_2'}^{\bar{l}_1}$. Then $\vph(y_1^{
k_1})={y_1'}^{ak_1}{y_2'}^{-\bar{l}_2k_1}$ and
$\vph(y_1^{l_1}y_2^{l_2})={y_1'}^l$.

In summary, by Corollary \ref{isorem}, we may assume that $l_2>0$
and $0\se l_1<l_2$. If $l_1=0$, then we go back to the case we just
discussed. If $l_1\neq0$ and $0<|q_{_1}|<1$, then $A$ is isomorphic to an
algebra in (V). Otherwise, $(\mc{D}, 0)$ is isomorphic to the datum
given by
\begin{itemize}
\item $\bgm=\lan y_1,y_2\ran\cong \ZZ^2$;
\item $g_1'=y_1^l, g_2'=y_1^{ak_1}y_2^{k_1\bar{l}_2}$.  $\bar{l}_1,\bar{l}_2\in \ZZ^{+}$  are the integers such that $l\bar{l}_1=l_1$, $l\bar{l}_2=l_2$, and $(\bar{l}_1,\bar{l}_2)=1$.  $a,b\in\ZZ$ are the integers  such that $a\bar{l}_1+b\bar{l}_2=1$ and $0<a<\bar{l}_2$.
\item  $\chi_1'(y_1)=q_{_1}^{-\frac{k_1\bar{l}_2}{l_2}}$,
$\chi_1'(y_2)=q_{_1}^{\frac{ak_1-l}{l_2}}$,
 $\chi_2'=\chi_1'^{-1}$.
\item $(a_{ij})$ is the Cartan matrix of type
$A_1\times A_1$;
\item $\lmd_{12}=0$
\end{itemize}
via the isomorphism $(\vph,(12),\al_1=\al_2=1)$, where $\vph$ is the
isomorphism defined by $\vph(y_1)=y_1^ay_2^{\bar{l}_2}$ and
$\vph(y_2)=y_1^by_2^{-\bar{l}_1}$. It follows that $A$ is isomorphic
to an algebra in (V) as well.

It is clear that the algebras from different cases and subcases are
non-isomorphic to each other. It is sufficient to show that the
algebras in the same subcases in Case 2 are non-isomorphic. Each
algebra in (II), (IV) and (VI) is a lifting of an algebra in (I),
(III) and (V) respectively. So it is sufficient to discuss the cases
(I), (III) and (V).

First we discuss the case (I). Let  $\mc{D}$ and $\mc{D}'$ be two
data given by
\begin{itemize}
\item $\bgm=\lan y_1,y_2 \ran\cong \ZZ^2$;
\item $g_1=g_2=y_1^{k}$ for some $k\in \ZZ^+$;
\item  $\chi_1(y_1)=q_{_1}$, $\chi_1(y_2)=q_{_2}$, where $q_{_1},q_{_2}\in\kk$
satisfy that $0<|q_{_1}| <1$ and $q_{_1}$ is not a root of unity, and
$\chi_2=\chi_1^{-1}$;
\item $(a_{ij})_{1\se i,j\se 2}$ is the Cartan matrix of type $A_1\times A_1$
\end{itemize}
and

\begin{itemize}
\item $\bgm=\lan y_1',y_2' \ran\cong \ZZ^2$;
\item $g_1=g_2=y_1^{k'}$ for some $k'\in \ZZ^+$;
\item  $\chi_1(y_1)=q_{_1}'$, $\chi_1(y_2)=q_{_2}'$, where
$q_{_1}',q_{_2}'\in\kk$ satisfy that $0<|q_{_1}'| <1$ and $q_{_1}'$ is not a
root of unity, and $\chi_2=\chi_1^{-1}$;
\item $(a_{ij}')_{1\se i,j\se 2}$ is the Cartan matrix of type $A_1\times A_1$
\end{itemize}
respectively.  Assume that $(\vph,\sg,\al)$ is an isomorphism from
$(\mc{D},0)$ to $(\mc{D}',0)$. Say $\vph(y_1)={y_1'}^a{y_2'}^c$ and
$\vph(y_2)={y_1'}^b{y_2'}^d$. Since $g_1=g_2$ and $g_1'=g_2'$, we
have $\vph(y_1^k)=y_1'^{k'}$. Moreover, $k,k'>0$. So $a=1$, $c=0$
and $d=\pm1$. Consequently, we have $k=k'$, $q_{_1}=q_{_1}'$. If
$\sg=\id$, then $q_{_2}'=q_{_1}^{-b}q_{_2}$. Otherwise,
$q_{_2}'=q_{_1}^{b}q_{_2}^{-1}$. We have identified the isomorphic algebras
in (I).

Similarly, it is direct to show that each triple $(k,l,q)\in
\ZZ^+\times\ZZ^+\times\kk$ such that $0<|q|<1$ determines a
non-isomorphic algebra in (III).

Now we show that the algebras in (V) are non-isomorphic. Let
$\mc{D}$ and $\mc{D}'$ be the data given by
\begin{itemize}\item $\bgm=\lan y_1,y_2\ran\cong\ZZ^2$;  \item $g_1=y_1^{k}$, $g_2=y_1^{l_1}y_2^{l_2}$
such that  $k,l_1,l_2\in \ZZ^+$ and   $0<l_1<l_2$;\item
$\chi_1(y_1)=q$, where $q\in\kk$  is not a root of unity, $0<|q|<1$,
and $\chi_1(y_2)=q^{\frac{k-l_1}{l_2}}$ and
$\chi_2=\chi_1^{-1}$;\item $(a_{ij})_{1\se i,j\se 2}$,   the Cartan
matrix of type $A_1\times A_1$\end{itemize} and
\begin{itemize}\item $\bgm'=\lan y'_1,y'_2\ran$ is  also a free
abelian group of rank 2;
\item $g'_1={y'_1}^{k'}$, $g'_2={y'_1}^{l'_1}{y'_2}^{l'_2}$ such
that $k',l'_1,l'_2\in \ZZ^+$  and   $0<l'_1<l'_2$;\item
$\chi'_1(y'_1)=q'$, where $q'\in\kk$  is not a root of unity,
$0<|q'|<1$, and $\chi'_1(y'_2)=q'^{\frac{k'-l'_1}{l'_2}}$ and
$\chi'_2={\chi'_1}^{-1}$;\item $(a_{ij})_{1\se i,j\se 2}$,  the
Cartan matrix of type $A_1\times A_1$\end{itemize}respectively. We
claim that  $(\mc{D},0)$ and $(\mc{D}',0)$ are isomorphic if and
only if $q=q'$, $k=k'$, $l_1=l'_1$ and $l_2=l_2'$.

Assume that $(\mc{D},0)$ is isomorphic to $(\mc{D}',0)$ via an
isomorphism $(\vph,\sigma,\al_1=\al_2=1)$.  Suppose that
$\vph(y_1)={y'_1}^a{y'_2}^c$ and $\vph(y_2)={y'_1}^b{y'_2}^d$, with
$a,b,c,d\in \ZZ$.

Either $\sg=\id$ or $\sg=(12)$. If $\sigma=\id$, then
$\vph(g_i)=g'_i$, $i=1,2$. So $${y'_1}^{ak}{y'_2}^{ck}={y'_1}^{k'}
\t{ and }
 {y'_1}^{al_1+bl_2}{y'_2}^{cl_1+dl_2}={y'_1}^{l'_1}{y'_2}^{l'_2}.$$
Since $\vph$ is an isomorphism, we have $ad-bc=\pm1$.  Because,
$k,k',l_2,l'_2>0$,  $0< l_1<l_2$ and $0< l'_1<l'_2$, it follows that
$b=c=0$ and $a=d=1$. Therefore,  $k=k'$,  $l_1=l'_1$, $l_2=l_2'$,
and $q=q'$. Namely, $(\mc{D},0)=(\mc{D}',0)$

If $\sigma=(12)$, then $\vph(g_i)=g'_{3-i}$, $i=1,2$. This implies
that $${y'_1}^{ak}{y'_2}^{ck}={y'_1}^{l'_1}{y'_2}^{l'_2} \t{ and }
 {y'_1}^{al_1+bl_2}{y'_2}^{cl_1+dl_2}={y'_1}^{k'}.$$ We can find
 integers $\bar{l}_1$ and $\bar{l}_2$, such that
$l_1=\bar{l}_1l$, $l_2=\bar{l}_2l$,  $l>0$ and
$(\bar{l}_1,\bar{l}_2)=1$.

Since $ad-bc=\pm1$, we have $(c,d)=1$. From $ck=l'_2>0$ and
$cl_1+dl_2=0$, it follows that $c=\bar{l}_2$ and $d=-\bar{l}_1$. If
$ad-bc=1$,  we have
$$k'=al_1+bl_2=l(a\bar{l}_1+b\bar{l}_2)=-l(ad-bc)=-l<0,$$ a
contradiction!

If $ad-bc=-1$, we have
$$ q'= \chi'_1(y_1')=\chi_2\vph^{-1}(y_1)=\chi_2(y_1^{\bar{l}_1}y_2^{\bar{l}_2})=q^{-\bar{l}_2\frac{k}{l_2}}.
$$ But $\bar{l}_2,k,l_2>0$ and
$0<|q|,|q'|<1$. We get a contraction as well. In summary, we have proved
the claim. \qed

Now we list all pointed CY Hopf algebras $U(\mc{D},\lmd)$  of
dimension 4  in terms of  generators and relations in the following
table. Note that $q_{_1}$ and $q$ are not roots of unity.
\begin{center}
CY algebras of dimension 4
\begin{longtable}{ccc}

\hline Case &Generators&Relations \\ \hline

Case 1& $y_h,y_h^{-1}$ &
$y_h^{\pm1}y_m^{\pm1}=y_m^{\pm1}y_h^{\pm1}$\\
&$1\se h\se 4$&$y_h^{\pm1}y_h^{\mp1}=1$\\
&&$1\se h,m\se 4$\\
\hline

Case 2 (I)&$y_1^{\pm1},y_2^{\pm1},x_1,x_2$&$y_h^{\pm1}y_m^{\pm1}=y_m^{\pm1}y_h^{\pm1}$\\
&&$y_h^{\pm1}y_h^{\mp1}=1$\\
&&$1\se h,m\se 2$\\
&&$y_1x_1=q_{_1}x_1y_1$, $y_1x_2=q_{_1}^{-1}x_2y_1$\\
&&$y_2x_1=q_{_2}x_1y_2$, $y_2x_2=q_{_2}^{-1}x_2y_2$\\&& $0<|q_{_1}|<1$\\
&&$x_1x_2-q_{_1}^{-k}x_2x_1=0$, $k\in \ZZ^+$\\
\hline

Case 2 (II)&$y_1^{\pm1},y_2^{\pm1},x_1,x_2$&$y_h^{\pm1}y_m^{\pm1}=y_m^{\pm1}y_h^{\pm1}$\\
&&$y_h^{\pm1}y_h^{\mp1}=1$\\
&&$1\se h,m\se 2$\\
&&$y_1x_1=q_{_1}x_1y_1$, $y_1x_2=q_{_1}^{-1}x_2y_1$\\
&&$y_2x_1=q_{_2}x_1y_2$, $y_2x_2=q_{_2}^{-1}x_2y_2$\\&& $0<|q_{_1}|<1$\\
&&$x_1x_2-q_{_1}^{-k}x_2x_1=1-y_1^{2k}$, $k\in \ZZ^+$\\
\hline

Case 2 (III)&$y_1^{\pm1},y_2^{\pm1},x_1,x_2$&$y_h^{\pm1}y_m^{\pm1}=y_m^{\pm1}y_h^{\pm1}$\\
&&$y_h^{\pm1}y_h^{\mp1}=1$\\
&&$1\se h,m\se 2$\\
&&$y_1x_1=qx_1y_1$, $y_1x_2=q^{-1}x_2y_1$\\
&&$y_2x_1=q^{\frac{k }{l }}x_1y_2$, $y_2x_2=q^{-\frac{k }{l }}x_2y_2$\\
&&$x_1x_2-q^{-k}x_2x_1=0$\\
&& $k,l \in \ZZ^+$, $0<|q|<1$\\
\hline

Case 2 (IV)&$y_1^{\pm1},y_2^{\pm1},x_1,x_2$&$y_h^{\pm1}y_m^{\pm1}=y_m^{\pm1}y_h^{\pm1}$\\
&&$y_h^{\pm1}y_h^{\mp1}=1$\\
&&$1\se h,m\se 2$\\
&&$y_1x_1=qx_1y_1$, $y_1x_2=q^{-1}x_2y_1$\\
&&$y_2x_1=q^{\frac{k }{l }}x_1y_2$, $y_2x_2=q^{-\frac{k }{l }}x_2y_2$\\
&&$x_1x_2-q^{-k}x_2x_1=1-y_1^{k }y_2^{l }$\\
&& $k,l \in \ZZ^+$, $0<|q|<1$\\
\hline

Case 2 (V)&$y_1^{\pm1},y_2^{\pm1},x_1,x_2$&$y_h^{\pm1}y_m^{\pm1}=y_m^{\pm1}y_h^{\pm1}$\\
&&$y_h^{\pm1}y_h^{\mp1}=1$\\
&&$1\se h,m\se 2$\\
&&$y_1x_1=qx_1y_1$, $y_1x_2=q^{-1}x_2y_1$\\
&&$y_2x_1=q^{\frac{k-l_1}{l_2}}x_1y_2$, $y_2x_2=q^{-\frac{k-l_1}{l_2}}x_2y_2$\\
&&$x_1x_2-q^{-k}x_2x_1=0$\\
&& $k,l_1,l_2\in \ZZ^+$, $0< l_1<l_2$, $0<|q|<1$\\
\hline

Case 2 (VI)&$y_1^{\pm1},y_2^{\pm1},x_1,x_2$&$y_h^{\pm1}y_m^{\pm1}=y_m^{\pm1}y_h^{\pm1}$\\
&&$y_h^{\pm1}y_h^{\mp1}=1$\\
&&$1\se h,m\se 2$\\
&&$y_1x_1=qx_1y_1$, $y_1x_2=q^{-1}x_2y_1$\\
&&$y_2x_1=q^{\frac{k-l_1}{l_2}}x_1y_2$, $y_2x_2=q^{-\frac{k-l_1}{l_2}}x_2y_2$\\
&&$x_1x_2-q^{-k}x_2x_1=1-y_1^{k+l_1}y_2^{l_2}$\\
&& $k,l_1,l_2\in \ZZ^+$, $0< l_1<l_2$, $0<|q|<1$\\
\hline

\end{longtable}
\end{center}

Let $\mathfrak{g}$ be a semisimple Lie algebra and
$U_q(\mathfrak{g})$ its quantized enveloping algebra.  By
\cite[Prop. 6.4]{bz}, the global dimension of the algebra
$U_q(\mathfrak{g})$ is the dimension of $\mathfrak{g}$. Thus, if
$U_q(\mathfrak{g})$ is of global dimension less than 5, then
$U_q(\mathfrak{g})$ is isomorphic to $U_q(\mathfrak{sl}_2)$, which
is of global dimension 3. That is, among the algebras of the form
$U_q(\mathfrak{g})$, only $U_q(\mathfrak{sl}_2)$ appears in the
lists of Propositions \ref{12}, \ref{3} and \ref{4}.  The algebra
$U_q(\mathfrak{sl}_2)$ is isomorphic to $U(\mc{D},\lmd)$ with  the
datum  given by
\begin{itemize}
\item $\bgm=\lan g\ran$, a free abelian group of rank 1;
\item The Cartan matrix is of type $A_1\times A_1$;
\item $g_1=g_2=g$;
\item $\chi_1(g)=q^{-2}$, $\chi_2(g)=q^{2}$, where $q$ is not a root of unity;
\item $\lmd_{12}=1$.\end{itemize}It belongs to (II) of Case 2 of
Proposition \ref{3}.

The family of pointed Hopf algebras $U(\mc{D},\lmd)$ provide more
examples of CY Hopf algebras of higher dimensions. From the
classification of CY pointed Hopf algebras $U(\mc{D},\lmd)$ of
dimensions less than 5, we see that the Cartan matrices are either
trivial or of type $A_1\times\cdots\times A_1$. The following
example provides a CY pointed Hopf algebra of type $A_2\times A_1$
of dimension  7.

\begin{eg}\label{eg pointed}
Let $A$ be $U(\mc{D},\lmd)$ with the datum $(\mc{D},\lmd)$   given
by
\begin{itemize}
\item $\bgm=\lan y_1,y_2,y_3 \ran$, a free abelian group of rank 3;
\item The Cartan matrix is $$\left(\begin{array}{ccc}2&-1&0
\\-1&2&0\\0&0&2\end{array}\right);$$
\item $g_i=y_i$, $1\se i\se 3$;
\item $\chi_i$, $1\se i\se 3$, are given by the following table, where $q$ is not a root of unity. $$\begin{tabular}{|c|c|c|c|}\hline
&$y_1$&$y_2$&$y_3$\\\hline $\chi_1$&$q$&$q^{-2}$&$q^4$\\\hline
$\chi_2$&$q$&$q$&$q^{-2}$\\\hline
$\chi_3$&$q^{-4}$&$q^2$&$q^{-4}$\\\hline

\end{tabular}$$
\item $\lmd=0$
\end{itemize}
In other words, $A$ is the  algebra with generators $x_i$,
$y_j^{\pm1}$, $1\se i,j\se 3$, subject to the relations
$$y_i^{\pm1}y_j^{\pm1}=y_j^{\pm1}y_i^{\pm1},\;\;\;y_j^{\pm1}y_j^{\mp1}=1,\;\;\; 1\se i,j\se 3,$$
$$y_j(x_i)=\chi_i(y_j)x_iy_j,\;\;\;1\se i,j\se 3,$$
$$x_1^2x_2-qx_1x_2x_1-q^2x_1x_2x_1+q^3x_2x_1^2=0,$$
$$x_2^2x_1-q^{-2}x_2x_1x_2-q^{-1}x_2x_1x_2+q^{-3}x_1x_2^2=0,$$
$$x_1x_3=x_3x_1.$$

The non-trivial liftings of $A$ are also CY.
\end{eg}

\vspace{5mm}

\subsection*{Acknowledgement}
We thank Dr. Jiwei He for his helpful comments and suggestions. This work forms part of the first author's PhD thesis at Hasselt University.
\vspace{5mm}

\bibliography{}

\end{document}